\documentclass[12pt]{article}

\usepackage{graphicx}
\usepackage{psfrag}
\usepackage{amsmath}
\usepackage{amssymb}
\usepackage[dvips]{color}
\newcommand{\oc}[1]{\overset{\circ}{#1}}

\def\ftf{facet-to-facet}
\def\tes{tessellation}

\def\pie{$\pi$-edge}
\def\pies{$\pi$-edges}

\def\la{\lambda}
\def\al{\alpha}
\def\i{\infty}
\def\R{\mathbb{R}}

\def\Pr{\mathsf P}

\def\VE{\mu_{VE}}

\def\VZ{\mu_{VZ}}
\def\VP{\mu_{VP}}
\def\VF{\mu_{VZ_2}}
\def\VR{\mu_{VZ_1}}
\def\VS{\mu_{VP_1}}

\def\EP{\mu_{EP}}

\def\EF{\mu_{EZ_2}}

\def\ZP{\mu_{ZP}}
\def\ZE{\mu_{ZE}}
\def\ZV{\mu_{ZV}}
\def\nuZR{\nu_1(Z)}
\def\nuZA{\nu_0(Z)}
\def\nuZF{\nu_2(Z)}

\def\nuFS{\nu_1(Z_2)}
\def\nuFC{\nu_0(Z_2)}

\def\PV{\mu_{PV}}

\def\nuPS{\nu_1(P)}
\def\nuPC{\nu_0(P)}

\def\XY{\mu_{XY}}

\def\cn{^{\hspace{0.07cm}\circ}}

\def\cc{^{\hspace{0.07cm}\circ\hspace{.08cm}\circ}}

\def\nc{^{\hspace{0.3cm}\circ}}

\def\ds{\displaystyle}
\def\tf{\tfrac}

\def\Lra{\Longrightarrow}
\def\Llra{\Longleftrightarrow}

\begin{document}
\begin{center}
\Large{\textbf{Constraints on the fundamental topological parameters of spatial tessellations}}\\
\small{\ }\\

Richard Cowan\\ School of Mathematics and Statistics\\ University
of Sydney, NSW, 2006, Australia\\ \textit{e-mail:
rcowan@usyd.edu.au}\\
\ \\
Viola Weiss\\ Ernst--Abbe--Fachhochschule Jena\\
D-07703 Jena, Germany.
\\ \textit{e-mail: Viola.Weiss@fh-jena.de
}
\end{center}
\date
\maketitle

\begin{abstract}
Tessellations of $\R^3$ that use convex polyhedral cells to fill the
space can be extremely complicated, especially if they are not
`\ftf', that is, if the facets of a cell do not necessarily coincide
with the facets of that cell's neighbours. In a recent paper
\cite{wc}, we have developed a theory which covers these complicated
cases, at least with respect to their combinatorial topology. The
theory required seven parameters, three of which suffice for \ftf\
cases; the remaining four parameters are needed for the awkward
adjacency concepts that arise in the general case. This current
paper establishes constraints that apply to these seven parameters
and so defines a permissible region within their seven-dimensional
space, a region which we discover is not bounded. Our constraints in
the relatively simple facet-to-facet case are also new.
\end{abstract}
\begin{flushleft}
\hspace{2.45cm}Secondary: 60G55; 51M20;  52B10
\end{flushleft}

\section*{1. Introduction}

In this paper we continue our study of \emph{random stationary spatial
tessellations}, that is, random \tes s of the three-dimensional
space $\R^3$ having statistical properties that are invariant under
translation. In \cite{wc}, we developed a theory of such \tes s in cases where the cells of the \tes\ are closed
 convex polyhedra and not necessarily facet-to-facet. Seven parameters were needed to address the major topological issues.

\textbf{Primitive elements:} Three of these parameters suffice for
the relatively simple facet-to-facet case, where the questions of
interest focus on the four primitive elements of the \tes,
\emph{vertices}, \emph{edges}, \emph{plates} and \emph{cells}. (We
use the word \emph{plate} for a closed convex polygon which lies on
the boundary of two cells; if these two closed cells are $z$ and
$z'$, then the polygon $z\cap z'$ is that plate.) The three
parameters, with cyclic subscripts in the letters $V, E$ and $P$,
are:
\begin{itemize}
    \item $\VE$, the expected number of edges emanating from the typical vertex;
    \item $\EP$, the expected number of plates emanating from a typical edge;
    \item $\PV$, the expected number of vertices lying on the boundary of the typical polygonal plate. Note that if a plate is an $n$-gon, this number may exceed $n$ in the non facet-to-facet case because plate \emph{corners} are not the same entities as \tes\ \emph{vertices}.
\end{itemize}

As a matter of notation, classes of the primitive elements
\emph{vertices}, \emph{edges}, \emph{plates} and \emph{cells} are
called $V, E, P$ and $Z$. The generic expression, $m_Y(x)$ for $x\in
X$, is defined as the number of objects of type $Y$ \emph{adjacent
to} a particular $x$ in the object class $X$. We define $\XY$ as the
expected number of $Y$-type objects \emph{adjacent to} the typical
object of class $X$. Two objects $x$ and $y$ are \emph{adjacent} if
either $x\subseteq y$ or $y\subseteq x$. Adjacency is a precise
concept which covers relationships like `emanating from' or `lying
on', used above. All twelve mean adjacencies $\XY$, for $X$ and $Y$
primitive, $X\ne Y$, can be expressed (see Table 1, proved in
\cite{wc}) as functions of the \emph{cyclic adjacency} parameters,
$\VE$, $\EP$ and $\PV$.

From \cite{wc}, we reproduce the following
table. It uses the abbreviation,
\begin{equation}\label{fx}
f(x) := \VE\ \EP- x\ (\VE-2).
\end{equation}
\begin{center}
\begin{tabular}{|c||c|}
  \hline
  $X$ & $\la_X/\la_V$ \\\hline\hline
  vertices $\overset{\ }{V}$ & 1 \\[.2cm]
  edges $E$ &  $\tfrac12\overset{\ }{\VE}$ \\[.3cm]
  plates $P$ & $\displaystyle{\frac{\overset{\ }{\VE}\ \EP}{2\PV}}$ \\[.2cm]
  cells $Z$ &  $\displaystyle{\frac{\overset{\ }{f}(\PV)}{2\PV}}$\\
  \hline
\end{tabular}
\hspace{1mm}
\begin{tabular}{|c||c|c|c|c|}
  \hline
  $\mu$    & $V$ & $E$ & $P$ & $Z$  \\ \hline \hline
  $\overset{\ }V$ &  1  & $\VE$  & $\tfrac12\VE\ \EP$ & $\tfrac12f(2)$   \\[.2cm]
  $E$ &  2  &   1   & $\EP$ & $\EP$\\[.3cm]
  $P$ & $\PV$&$\PV$&1&$\overset{\ }{2}$\\[.3cm]
  $Z$ & $\displaystyle{\frac{\overset{\ }{\PV}f(2)}
      {f(\PV)}}$
  & $\displaystyle{\frac{\VE\ \EP\ \PV}{f(\PV)}}$
  & $\displaystyle{\frac{2\overset{\ }{\VE}\ \EP}{f(\PV)}}$ &  1 \\
  \hline
\end{tabular}

\end{center}
\noindent {\small Table 1:} {\scriptsize  Primitive intensities and adjacencies,
using the
abbreviation $f$, defined above in (\ref{fx}). All entries in the table are
expressed in terms of the scale parameter, $\la_V$, and the cyclic three
mean-adjacencies. Note that there are three linear identities within the
table: $\la_V-\la_E+\la_P-\la_Z = 0$; $\VE -\VP+\VZ=2$ ; $\ZV-\ZE+\ZP =2$.
These have links with Euler's polyhedral formula and with related formulae
from cell-complex theory (\cite{lz}, \cite{wz}).}\\

The intensity of objects of type $X$ (mean number of centroids per
unit volume) is denoted by $\la_X$. The value of $\la_V$ determines
the scale of the \tes\ and in Table 1  the three others, $\la_E,
\la_P$ and $\la_Z$, are also expressed in terms of $\la_V$ and the
cyclic adjacencies, $\VE$, $\EP$ and $\PV$.

Subsets of the generic class $X$ are denoted by $X[\cdot]$, with a
suitably chosen symbol in the square brackets. For example, shortly
we introduce a subclass of edges known as \pies; we use the notation
$E[\pi]$ for this subclass.

Adjacency is a symmetric relationship and an important identity
applies, proved in M\o ller's Theorem 5.1 \cite{mol}  and discussed
as equation (5) in \cite{wc}.
\begin{equation}\label{sym}
\la_X \XY = \la_Y \mu_{YX}.
\end{equation}

Although the results in Table 1 have been known for some time, at
least since 1980  (see \cite{rad}), there have been no studies of
the constraints which apply to the three mean-adjacencies,  $\VE$,
$\EP$ and $\PV$ (even in the facet-to-facet situation). We shall
rectify  this deficiency in this paper. After presenting some
examples of spatial \tes s in Section 2, we then focus on \ftf\
cases  and show, in Section 3, that:
\begin{align}
        4\leq&\ \VE; \qquad
        3\leq\EP \leq 6\Bigl(1-\frac{2}{\VE}\Bigr)<6;\qquad
        3\leq\PV < \frac{\VE\EP}{\VE-2}\le 6.\label{case1}
\end{align}
So these constraints, which are derived using the information in
Table 1 combined with elementary geometry of convex polyhedra, apply
to \emph{all} facet-to-facet \tes s. In Section 4, we present more
\ftf\ examples, including one construction in which $\VE$ can be
arbitrary large. Thus we reinforce the absence of an upper bound for
$\VE$ in (\ref{case1}).

\begin{figure}[ht]
   \psfrag{EP}{\small$\EP$}\psfrag{VE}{\small$\VE$}
\psfrag{146}{\scriptsize \hspace{-4.5mm}$1$\&$4$\&$6$d}
\psfrag{359}{\scriptsize \hspace{-3mm}$3$\&$5$\&$9$c}
\psfrag{ad}{\scriptsize \hspace{-4mm}$9$a\&$9$d}
\psfrag{a}{\scriptsize a}
\psfrag{b}{\scriptsize b}
\psfrag{c}{\scriptsize c}
\psfrag{d}{\scriptsize d}
\psfrag{e}{\scriptsize e}
\psfrag{f}{\scriptsize f}
\psfrag{g}{\scriptsize g}
\psfrag{i}{\scriptsize i}
\psfrag{h}{\scriptsize h}
\psfrag{10}{\scriptsize\hspace{-1.2mm}$10$}
\psfrag{11}{\scriptsize\hspace{-1.2mm}$11$}
\psfrag{12}{\scriptsize$12$}
\psfrag{13}{\scriptsize\hspace{-1.2mm}$13$}
\psfrag{14}{\scriptsize\hspace{-1.2mm}$14$}
\psfrag{15}{\scriptsize$15$}
\psfrag{16}{\scriptsize\hspace{-1mm}$16$}
\psfrag{17}{\scriptsize$17$}  \psfrag{18}{\scriptsize\hspace{-1.2mm}$18$}
\psfrag{1}{\scriptsize\hspace{-.6mm}$1$}
\psfrag{2}{\scriptsize\hspace{-.6mm}$2$}
\psfrag{3}{\scriptsize\hspace{-.6mm}$3$}
\psfrag{4}{\scriptsize\hspace{-.6mm}$4$}
\psfrag{5}{\scriptsize\hspace{-.6mm}$5$}
\psfrag{6}{\scriptsize\hspace{-.6mm}$6$}
\psfrag{7}{\scriptsize\hspace{-.6mm}$7$}
\psfrag{8}{\scriptsize\hspace{-.6mm}$8$}
\psfrag{9}{\scriptsize\hspace{-.6mm}$9$}
    \begin{center}
        \includegraphics[width=163mm]{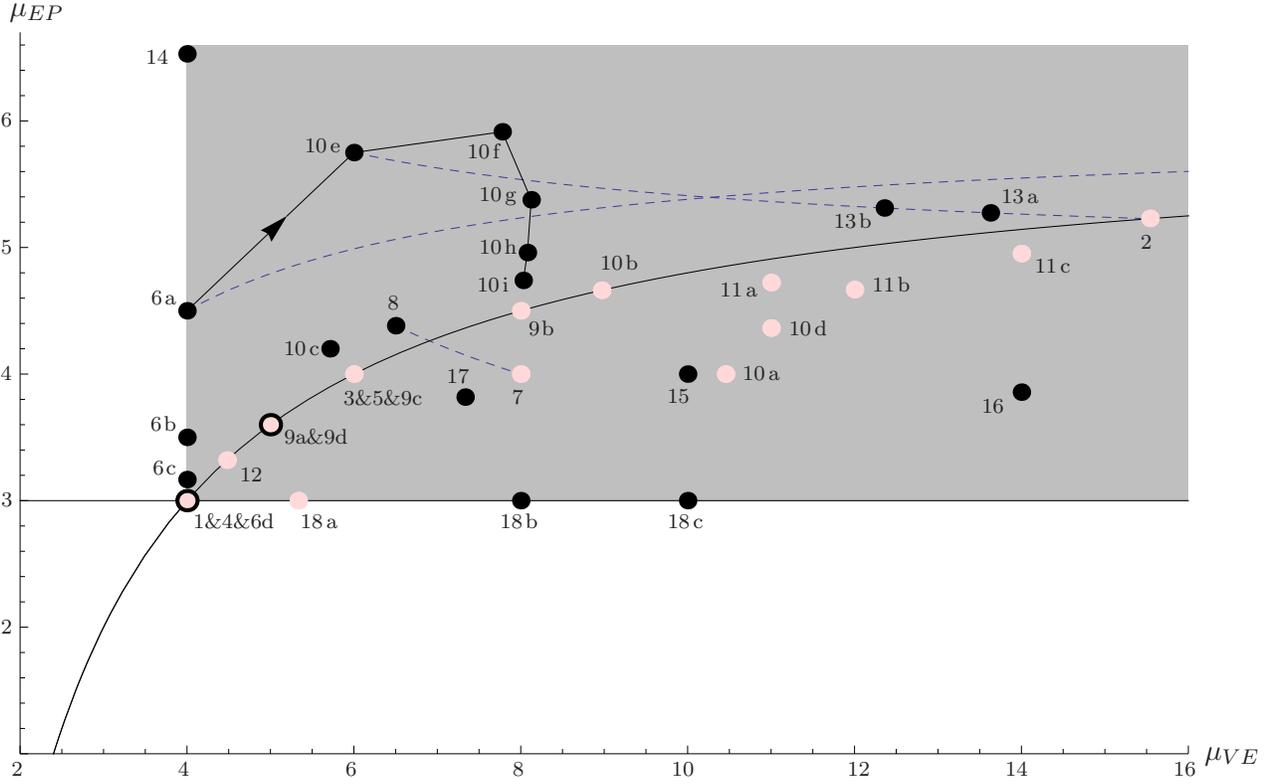}\\
        \caption{\scriptsize \label{fig1} A plot of
        $\EP=6(1-2/\VE)$, shown as the solid smooth curve, together
        with $(\VE, \EP)$ values for the various models discussed in this
        paper. Facet-to-facet models have a light-coloured dot and non \ftf\ models a
        black dot. Two of the dots are hybrid dots; the $(4,3)$ dot represents one
        \ftf\ model (Example 1) and two that are not (Examples 4 and 6d) and, for the
        $(5,18/5)$ dot, Example 9a is \ftf\ whilst 9d is not. Dashed curves and the
         connected sequence of line segments joining some of the dots are explained in the text.}
     \end{center}
\end{figure}

The examples throughout our paper enable us to overlay data points
$(\VE,\EP)$ on a plot of the curve $\EP= 6(1-2/\VE)$, revealing the
role that the middle constraint in (\ref{case1}) plays. We shall
discover, as illustrated in Figure \ref{fig1}, that one can find
some \emph{non} facet-to-facet examples where (\ref{case1}) also
holds --- though for many such examples (\ref{case1}) is violated.

The dashed curves in Figure \ref{fig1} result from \emph{mixtures of
\tes s}; these are explained in Section 5. The dots which are linked
by a connected sequence of line segments are explained in Example
10.

\textbf{Faces of the primitive elements:} Section 6 commences the
theory for \tes s which are not \ftf. It presents a different set of
constraints for $\VE, \EP$ and $\PV$ as well as constraints on the
four additional parameters, $\xi, \kappa, \psi$ and $\tau$. These
four parameters come from our earlier study \cite{wc} and quantify
the consequences of having cell-facets that do not coincide with
those of their neighbours.

To define these parameters, we need to settle the terminology for
faces of the primitive elements --- in view of the fact that the
standard terminologies in \tes\ and polytope theories clash. We have
decided to use the \tes\ theoretic meaning of words \emph{vertex},
\emph{edge}, \emph{plate} and \emph{cell}; other names, \emph{apex,
ridge} and \emph{facet} are used for the $j$-faces of a cell,
$j=0,1,2$. A j-face of a polygon (j=0,1) we call a \emph{corner} or
\emph{side}. We distinguish between sides of plates and of facets by
using the terms plate-side and facet-side.

We note that in a \ftf\ \tes\ every face of a primitive element is
equal (as a set in $\R^3$) to  a primitive element. For example,
every ridge ($1$-face of a cell) equals an edge of the \tes\ -- and
every edge equals at least three ridges. This simplicity is lost in
the non \ftf\ case.

Let $X$ be a class of convex polytopes, each member of the class
having dimension $i\leq 3$. We define $X_j, j<i$, as the class of
objects which are $j$-dimensional faces ($j$--faces) of some
polytope belonging to $X$. For instance $P_1$ is the class of the
\emph{plate-sides} (1-faces of plates) and $Z_2$ is the class of
\emph{facets} ($2$-faces of cells). We denote the expected number of
$j$-faces that a typical polytope of class $X$ has by $\nu_j(X)$.
For example, $\nuZR$ is the expected number of ridges ($1$-faces of
the cells) for a typical cell.  Adjacency formulae and intensities
involving these classes of elements were derived in our earlier
paper \cite{wc} and will be introduced in Section 6 as needed.
Section 2 of that paper also has a more detailed discussion of our
\emph{face--of--polytope} classes, which may in some cases be
\emph{multisets}.

\begin{figure}[ht]
    \begin{center}
    \psfrag{P}{\scriptsize $P$} \psfrag{C1}{\scriptsize\hspace{-2mm}$C_1$}
    \psfrag{C2}{\scriptsize \hspace{-1mm}$C_2$}
        \includegraphics[width=44mm]{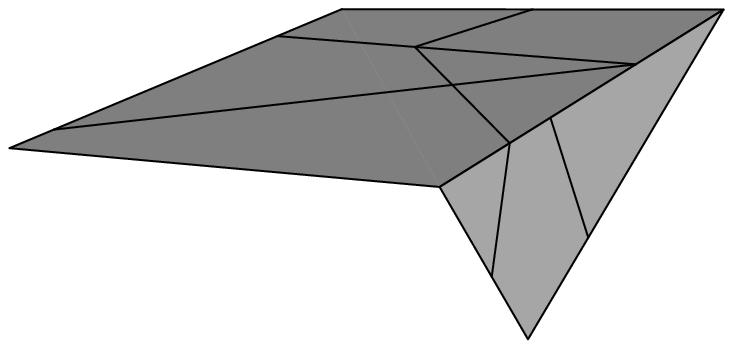}\hspace{.8cm}
        \includegraphics[width=58mm]{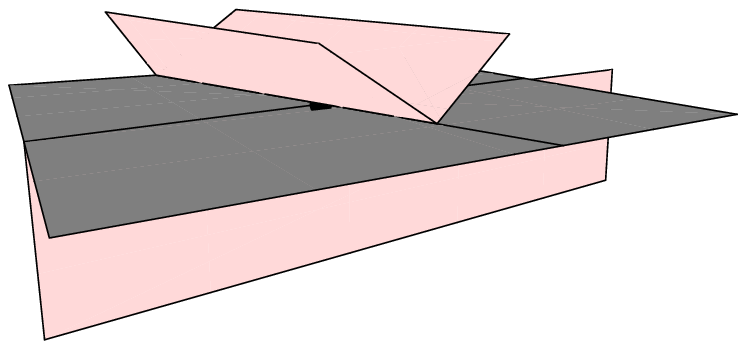}\hspace{-.4cm}
         \includegraphics[width=58mm]{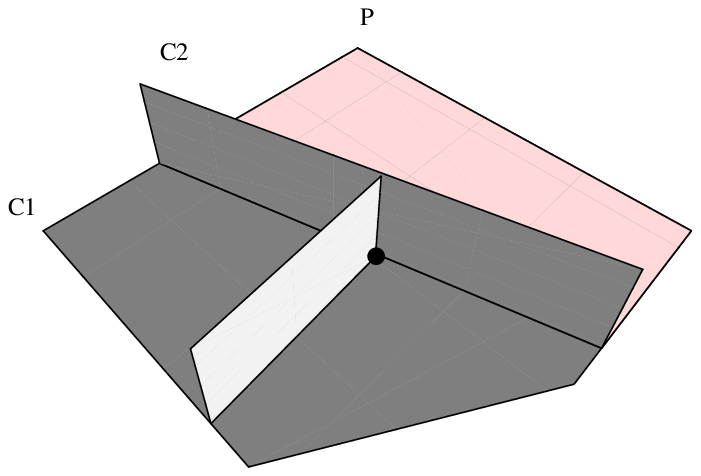}
     \vspace{-1cm}
          {\scriptsize \\(a) \hspace{3cm}(b) \ \ \includegraphics[width=15mm]{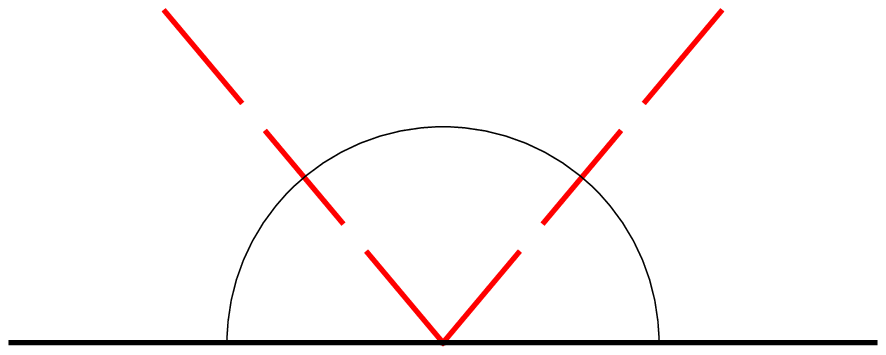} and \includegraphics[width=15mm]{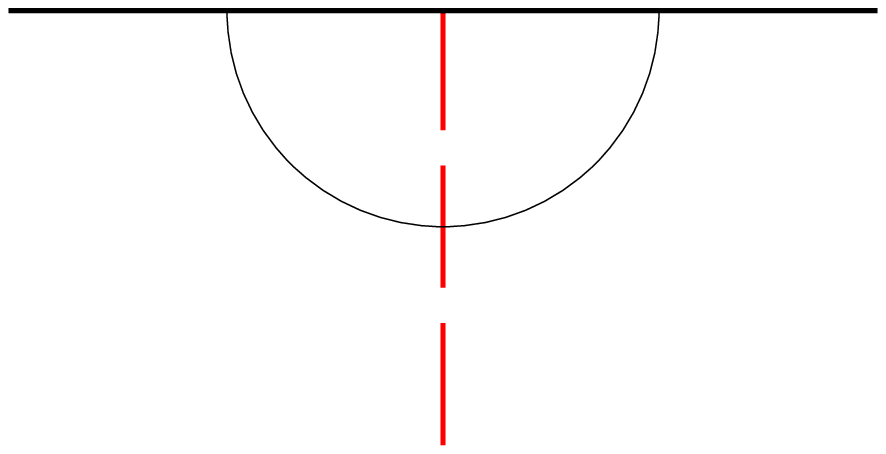}\hspace{2cm}(c)\ \ } \includegraphics[width=15mm]{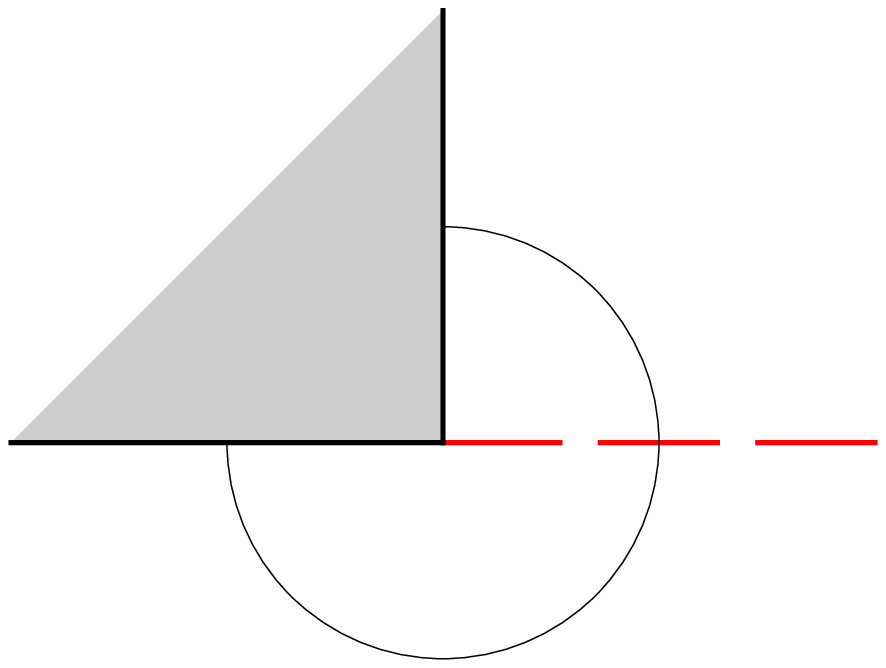}
        \caption{{\scriptsize \label{fig2} The schematic diagrams below the main figures are helpful in lemmas involving $\psi$ and $\tau$, as will be explained in Section 7 for the reader who does not guess their meaning immediately. The main figures are as follows. (a) Two facets
        of a cell $z$ are shown, with
        their interior structure due to neighbouring cells of $z$ (not shown).
        The dark-shaded facet (comprising six plates) has two
        hemi-vertices in its interior and seven $\pi$-edges. The lighter
        facet
        (comprising three plates) demonstrates that $\pi$-edges can
        exist on a facet without there being hemi-vertices. A ridge
        is also illustrated, with three interior vertices, two
        of which lie in the interior of a plate--side (whose plate is
        one of the nine
        shown). (b) A vertex adjacent to
        three plate--side-interiors and five ridge-interiors, but no
        facet-interior.  (c) A common type of hemi-vertex, with $P$ and $C_1$ being coplanar. The vertex is in the interior of one plate--side (that of the plate $P$) and one ridge (where $P$ and $C_2$ meet). The labels $C_1$ and $C_2$ are \emph{bookcovers}, a concept explained later. }}
    \end{center}
\end{figure}

\textbf{The four \emph{interior} parameters, $\xi, \kappa, \psi$ and
$\tau$:} The primitive elements, edge, plate and cell,  can have no
other element (or face of an element) lying in their \emph{relative
interior}. The face of a primitive element, however, can have such
interior structure --- and the four \emph{interior parameters}
quantify the prevalence of this phenomenon in non facet-to-facet
\tes s. These additional parameters introduced in \cite{wc} and
visualised in Figure \ref{fig2} (and especially in Figure 3 of
\cite{wc} in the context of elaborate cell architecture) are defined
as follows. We use the notation where $\oc{X}$ means the class of
\emph{relative interiors} of members of $X$, and where the word
`interior' will henceforth mean `relative interior'.
\begin{itemize}
    \item $\xi :=$ the proportion of edges in the \tes\ whose interior
    is contained in the
    interior of some facet --- these being called \emph{$\pi$-edges} in
    \cite{wc}. Thus $\xi = \la_{E[\pi]}/\la_E$. Note also that  an adjacency representation exists, namely $\xi = \EF\cc$.
    \item $\kappa :=$ the proportion of vertices in the \tes\
    contained in the interior of some facet --- these being
    called \emph{hemi-vertices}  because half of the
    neighbourhood of such a vertex is a hemi-sphere lying
    within one cell. Thus $\kappa = \la_{V[hemi]}/\la_V$. Note also that $\kappa = \VF\nc$.
    \item $\psi := \VR\nc$, the expected number of ridge-interiors adjacent to a typical vertex.
    \item $\tau := \VS\nc$, the expected number of plate--side-interiors adjacent to a typical vertex.
\end{itemize}
Being proportions, $\xi$ and $\kappa$ must lie in the unit interval, so
$(\kappa, \xi) \in [0,1]^2$. We note that the first of these four
parameters, $\xi$, is a descriptor of edges, whilst the last three describe
vertices. It is possible however to write $\xi$ in terms of the vertex
adjacency, $\mu_{VE[\pi]}$. Using (\ref{sym}), Table 1 and the fact that
$\mu_{E[\pi]V}=2$,
\begin{equation}\label{xipi}
    \xi := \frac{\la_{E[\pi]}}{\la_E} = \frac{\la_{E[\pi]}\,\mu_{E[\pi]V}}{2\la_E}=\frac{\la_V\,\mu_{VE[\pi]}}{2\la_E}=\frac{\mu_{VE[\pi]}}{\VE}.
\end{equation}
Later we use this linkage between the vertex property
$\mu_{VE[\pi]}$, the mean number of \pies\ adjacent to a typical
vertex, and the proportion  $\xi$ of \pies.

\textbf{The relationship with cell complex theory:} In the theory of
cell complexes  the terminology \emph{face-to-face} has been used
(see \cite{schw}, p.447). In the context of spatial \tes s, this
means the following:  if the intersection of two cells has dimension
$j\in \{0,1,2\}$, then the intersection is a $j$-face of both cells.
We have chosen to use the concept \emph{\ftf}, which has the same
meaning but apparently only when $j=2$. The following lemma provides
an equivalence of these two concepts.

{\sc Lemma 1:} \emph{Our concept of facet-to-facet is equivalent to
the concept of face-to-face used in the literature of cell
complexes.}

{\sc Proof:} Clearly, \emph{face-to-face} implies \emph{\ftf}, so we
focus on the converse to this.

Firstly, being \emph{facet-to-facet} is equivalent to having
$\xi=0$, because $\xi>0 \Llra$ \emph{the existence of \pies}\ and
having \pies \ would violate the condition that every facet of a
cell $z$ coincides with the facet of one of its neighbouring cells,
for all $z\in Z$.

Secondly, $\xi =0$ implies that $\kappa = \psi = \tau =0$. This
follows because hemi--vertices, and vertices which lie in the
interior of a ridge or a plate-side, have emanating \pies. If there
are no \pies, there can be none of these vertices.

Thirdly,  the condition $\xi= \kappa = \psi = \tau =0$  is
equivalent to being \emph{face-to-face}. This equivalence can be
understood by noting that the face-to-face condition prevents a
$j$-face of a cell $z$ being contained in the interior of a $k$-face
($0\leq j<k \leq 2$) of any other cell $z'$ (because $z\cap z'$
would be of dimension $j$ but not a $j$-face of $z'$.) Furthermore,
every vertex of the \tes\ coincides with some cell apices (cellular
$0$-faces) and every \tes\ edge is a subset of some cell ridges
(cellular $1$-faces). So a vertex cannot lie in a facet- or
ridge-interior and an edge-interior cannot lie in a facet-interior.
We further note that a plate-side (an entity which is not always a
face of a cell) is always contained in a ridge (which is), so also a
vertex cannot lie in a plate-side interior. So the lemma is now
proven. \hfill$\square$

Non \ftf\ \tes s are not cell complexes and that perhaps explains
why their study has been neglected.

\textbf{The general theory and our findings:} Section 6 shows how
the \emph{interior parameters} are constrained for those non
facet-to-facet \tes s\ conforming with (\ref{case1}). There are non
\ftf\ cases which do not conform to (\ref{case1}), because the upper
bound for $\EP$ no longer holds and the lower bound for $\PV$
becomes more complicated. Sections 7 and 8 also deals with these
cases; indeed one general theory presented in these sections covers
all the non \ftf\ cases.

It is convenient now to inform the reader of the major findings of
our general theory.

\vspace{.5cm} {\sc Theorem 1: Main result.}\\
\indent\emph{\underline{Facet-to-facet case}: For \tes s which are
\ftf, the constraints on the three cyclic parameters $\VE, \EP$ and
$\PV$ are as follows:}
\begin{align}
        4\leq&\ \VE; \qquad
        3\leq\EP \leq 6\Bigl(1-\frac{2}{\VE}\Bigr)<6;\qquad
        3\leq\PV < \frac{\VE\EP}{\VE-2}\le 6.\label{case1th}
\end{align}
\emph{In this \ftf\ case, the interior parameters $\xi, \kappa, \psi$ and $\tau$ all equal zero. }

\emph{\underline{Non \ftf\ case}: A \tes\ is not \ftf\ if and only
if $\xi>0$.}

\emph{Some tessellations which are not \ftf, have cyclic parameters
which conform to (\ref{case1th}), but many do not. The precise
constraints in this case are as follows:}
\begin{equation}\label{case2}
        4\leq\ \VE; \quad
        3\leq\EP;\quad
\begin{cases}
\qquad\quad\quad\ \ 3\leq \PV <
\ds{\frac{\VE\EP}{\VE-2}}<6& \mathrm{\ \ if\ }\EP<6\Big(1-\ds{\frac{2}{\VE}}\Bigr);\\
\ds{\frac{\VE\EP}{2(\VE-2)}}< \PV <  \ds{\frac{\VE\EP}{\VE-2}}&
\mathrm{\ \ if\ }\EP \geq 6\Big(1-\ds{\frac{\overset{\
}{2}}{\VE}}\Bigr).
\end{cases}
\end{equation}
\emph{The lower bound for $\PV$ in the second case is $> 3$ if and
only if $\EP > 6(1-2/\VE)$ whilst the upper bound for $\PV$ in this
case is $>6$.}

\emph{The permitted ranges for $\psi$ and $\tau$, given $\VE, \EP$
and $\PV$, are as follows.}
\begin{equation}\label{psibound}
    0\leq \psi \leq
\begin{cases}
\ \VE-2+\displaystyle{\frac{\VE\EP}{2}\Bigl(1-\frac{4}{\PV}\Bigr)}&\qquad \mathrm{if\ }\PV \leq \displaystyle{\frac{2\VE\EP}{3\VE-8}} \\\vspace{2mm}
    \ \displaystyle{\frac{\VE}{4}+\frac{\VE\EP}{2}\Bigl(1-\frac{3}{\PV}\Bigr)}&\qquad \mathrm{if\ }\PV \geq \displaystyle{\frac{\overset{\ }{2}\VE\EP}{3\VE-8}}.
\end{cases}
\end{equation}
\begin{align}
    \max\Big[0, \psi-\frac{\VE}{2}, \frac{\psi}{2} +\frac{\VE}{4}\Bigl(\EP-6\bigl(1-\frac{2}
    {\VE}\bigr)\Bigr)\Big] &\leq \tau \leq \min\bigl[\psi,\tfrac12
  \VE\EP
  \Bigl(1-\frac{3}{\PV}\Bigr)\bigr].\label{taubound}
\end{align}
\emph{The permitted range for $\kappa$ and $\xi$,
given $\VE, \EP, \PV, \psi$ and $\tau$, is given by}
\begin{align*}
        0\leq \kappa \leq\ \min\Bigl[1,\VE-2+\tfrac12
  \VE\EP
  \Bigl(1-\frac{4}{\PV}\Bigr)&-\psi\Bigr]
  \qquad\mathrm{and}\\
    \max\Bigl[\frac{2(\psi-\tau)+3\kappa}{\VE},\frac{4\psi+ 6\kappa}{\VE}-2
    \EP\Bigl(1-\frac{3}{\PV}\Bigr)\Bigr] \leq&\ \xi\ \leq
    \min\Bigl[1,
    \ 3 -\frac{\EP}{2}+ \frac{\psi - 6+3\kappa}{\VE}\Bigr],
\end{align*}
\emph{supplemented by $\xi>0$.}

\begin{figure}[ht]
\psfrag{tau}{$\tau$} \psfrag{psi}{\hspace{-1mm}$\psi$}
\psfrag{xi}{$\xi$} \psfrag{kap}{\hspace{-1mm}$\kappa$}
    \psfrag{EP}{{ \scriptsize $\hspace{-4mm}\EP$}}
        \psfrag{VE}{{ \scriptsize $\hspace{-4mm}\VE$}}
    \psfrag{PV}{{ \scriptsize$\hspace{-3.6mm}\PV$}}
    \begin{center}
        \includegraphics[width=165mm]{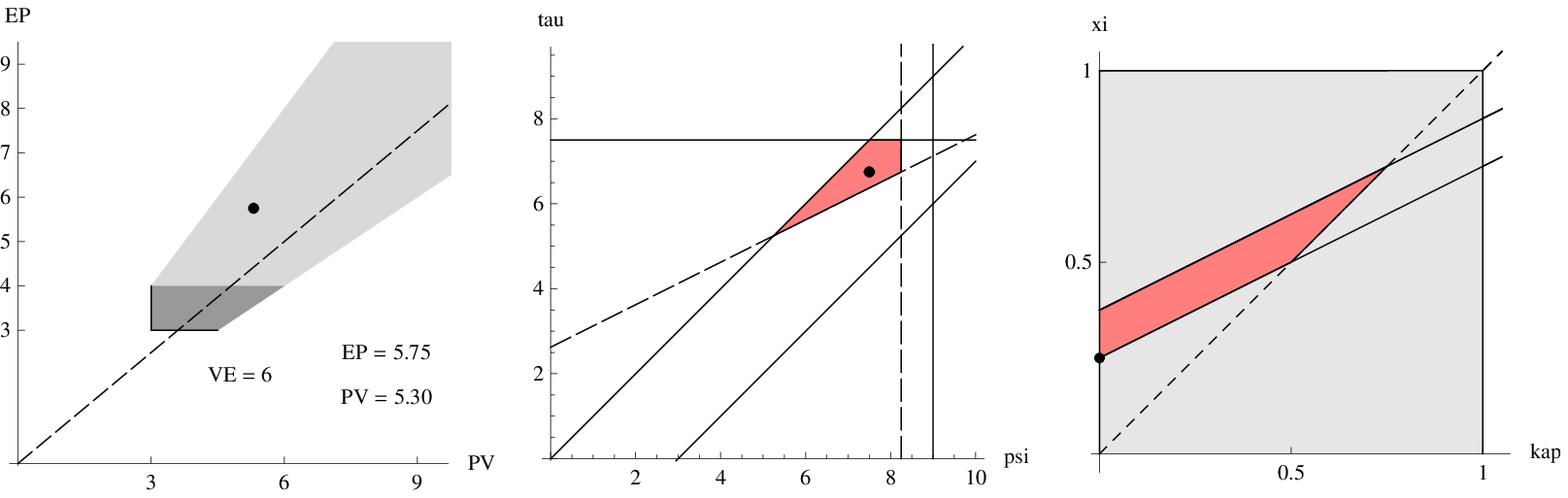}\\
          {\scriptsize (a)\hspace{5cm}(b)\hspace{5cm}(c)}
        \caption{{\scriptsize     \label{fig3} (a) The permitted range for
        $(\PV,\EP)$ if $\VE =6$. The dark shading, bounded above by the line
        $\EP = 6(1 -2/\VE)$, is the range for \ftf\ \tes s. It is `open' on
        its right boundary, but `closed' on the other boundaries.
        Some \tes s which are not \ftf\ have $(\PV,\EP)$ in the dark zone
        but, for others, the point lies in the light grey zone which is an
        open set, unbounded above.  The dashed line divides the region into
        two parts based on the inequalities in (\ref{psibound}). The black
        dots in the three diagrams correspond to the actual values of
        Example 10(e), described later. (b) Given $\VE=6$ and the $(\PV,\EP)$
        dot in (a), the dark region shows the permitted region for $(\psi,\tau)$. The constraints, which are straight lines, come from (\ref{psibound})  and (\ref{taubound}).  (c) The point $(\kappa, \xi)$ lies in the light grey square $[0,1]\times (0,1]$, but further contraints given in Theorem 1 lead to a smaller region, shown in darker shading. This plot assumes that the values of the other five parameters have been given.}}
    \end{center}
\end{figure}

\textbf{Our strategy:} It is not easy to present graphically the
shape of a seven-dimensional constraint space. Our method to give
some visual experience of the domain displays $\EP, \PV$ and the
four interior parameters in three two-dimensional plots: the $\EP$
versus $\PV$ permitted region given $\VE$; the $\tau$ versus $\psi$
permitted region, once the three cyclic parameters are given; the
permitted $\kappa$ versus $\xi$ region once the other five
parameters have been decided. This follows the way Theorem 1 has
been organised. A plot based on our Example 10(e), is presented as
Figure \ref{fig3}.

For the $\tau$ versus $\psi$ region, the boundaries  are determined firstly
by the cyclic parameters and secondly by the requirement that the region
does not make the $\xi$ versus $\kappa$ region null. Likewise, the region
$\EP$ versus  $\PV$ depends firstly on $\VE$ and secondly on the need to
have non-null regions for the interior parameters. As for $\VE$ itself,
we also investigate if there is any  value it can take in its range $[4,\i)$ that renders null the regions for the other six parameters. The reader will see these strategies used in the arguments that lead to Theorem 1.

Note that whilst we prove Theorem 1, space does not permit us to
prove all of the facts quoted in our example base. Separate papers
or supportive technical notes (\cite{cw1}, \cite{nwc}) provide some
proofs in the most substantial examples; they also give further
extension and graphical demonstrations of these examples.

{\sc Remark 1:} \textit{Before commencing the agenda set out above,
we note some very basic information. In a \tes\ with convex cells,
the cells must be polyhedra. Also, every vertex of the \tes\ must
have at least four emanating edges and every edge at least three
emanating plates; otherwise, non-convex cells would exist. So,
$\VE\geq 4$ and $\EP \geq 3$. Plates are polygons and therefore have
no fewer than three corners, so we can state that $\PV\geq 3$. Thus
the left-hand bound of the inequalities in (\ref{case1th}) and, in
part, (\ref{case2}) are proved trivially.}

{\sc Remark 2:} \textit{Cells are convex polyhedra, so we can use
some simple inequalities that apply to polyhedra. For
any convex polyhedron with $f_k$ faces of dimension $k$, and
therefore for any cell of the spatial \tes,}
\begin{align}
    3f_0\ \leq&\ 2f_1\label{simple}\quad \mathrm{and\ }\\
    3f_2\ \leq&\ 2f_1\label{simplex}.
\end{align}

\section*{2. Examples}

\begin{figure}[ht]
        \psfrag{A}{\hspace{-1mm}$A$}
        \psfrag{B}{\hspace{-.5mm}$B$}
    \begin{center}
    \includegraphics[width=50mm]{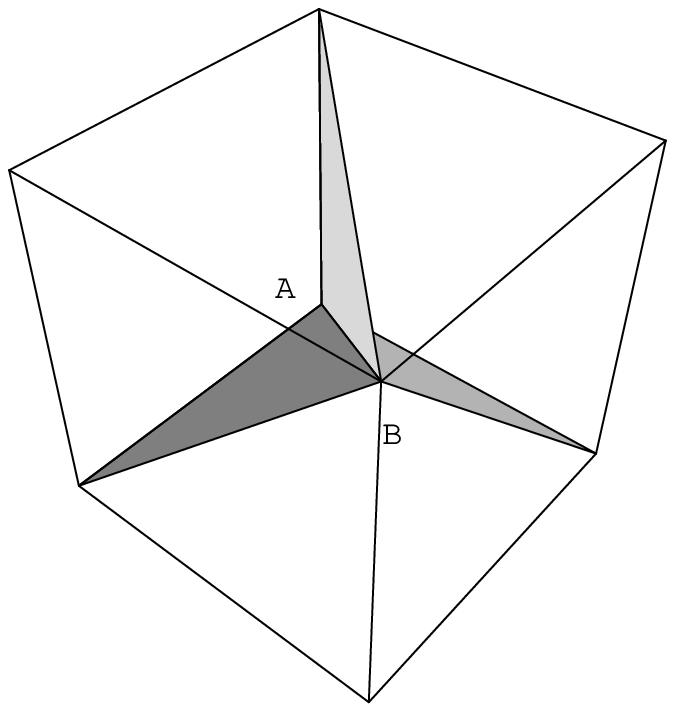}
    \includegraphics[width=50mm]{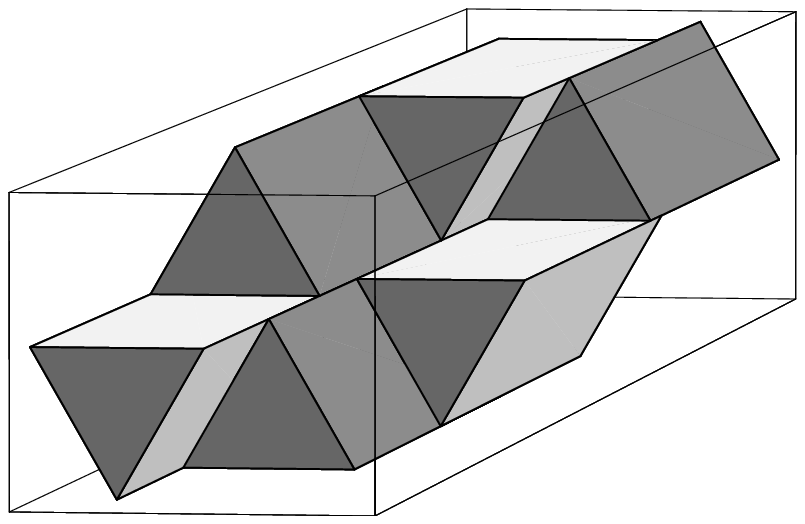} \hspace{.5cm}
    \includegraphics[width=50mm]{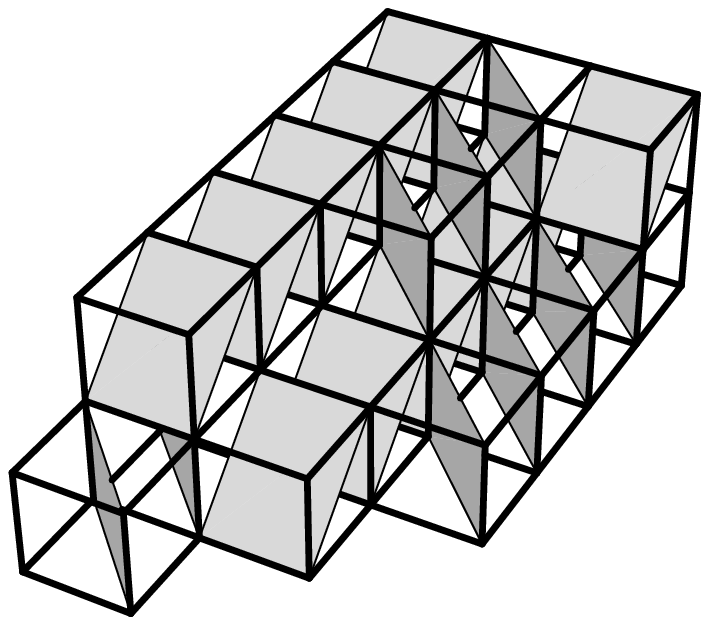} \\ \vspace{.3cm}
                 \renewcommand{\baselinestretch}{.1}\small
          {\scriptsize (a)\hspace{5cm}(b) \includegraphics[width=15mm]{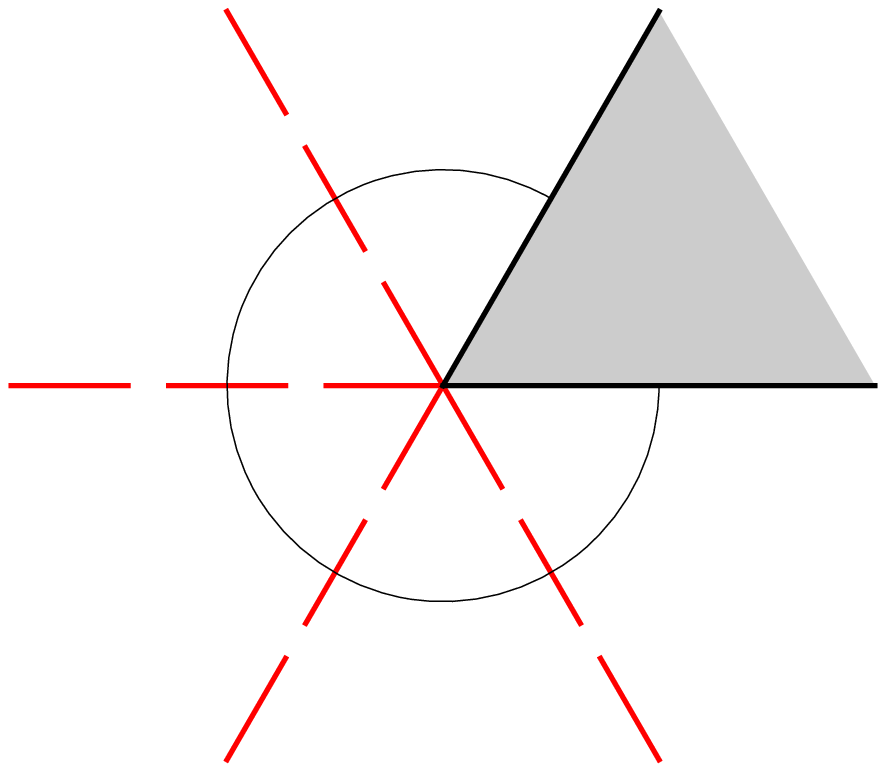}\hspace{5cm}(c)}
        \caption{\scriptsize \label{fig4} An illustration of
        (a) a partition of the cube used in Examples 7 and 16,
        (b) Example 6a (with the schematic, which applies to the visually
        most left vertex of the most right prism, defined in Section 7)  and
        (c) Example 15.}
    \end{center}
\end{figure}

We commence with well-known examples of random stationary \tes s in $\R^3$.
\begin{itemize}
    \item\textbf{ Example 1:} The Voronoi \tes\ based on seeds
        from a stationary Poisson point process. Here $(\VE,\EP) = (4,3)$ and $\PV=144\pi^2/(24\pi^2+35)\approx 5.2$. This is the best known \ftf\ \tes\
        (see \cite{okabe}).
    \item\textbf{ Example 2:} The Delaunay \tes, with similar seeds (also see \cite{okabe}). Here we have $(\VE, \EP)=(2+48\pi^2/35, 144\pi^2/(24\pi^2+35))\approx (15.5, 5.2)$. All plates are triangular and all cells are tetrahedra (so, being a \ftf\ \tes, $\PV=3$ whilst $\ZV = \ZP =4$ and $\ZE=6$).
    \item\textbf{ Example 3:} The \tes\ formed by random planes, no four of which meet at a point (\cite{mil}). Note that $(\VE,\EP) = (6,4)$. Also $\PV=4$ if the plane process is Poisson.
    \item\textbf{ Example 4:} The STIT \tes\ in $\R^3$ (see \cite{nw1}). This is our first example that is not \ftf. It is known from \cite{nw1} that $\VE=4, \EP=3$ and $\PV=\tf{36}{7}$ and from \cite{wc} (also \cite{thw}) that $\xi=1, \kappa=\tf23, \psi = 2$ and $\tau=\tf43$.
    \item \textbf{ Example 5:} Cubes packed in a lattice. Note that all
    tilings based on a repeating sub-unit can be converted to a random
    stationary \tes\ by locating the origin uniformly distributed within the
    sub-unit. The cyclic parameters are the same as in Example 3.
\end{itemize}
Each of the Examples 1--5 provide dots for the plot of Figure
\ref{fig1}; they are annotated by `example number'. Also on the
figure are dots which belong to a \emph{class of models}. For
example, we shall introduce our sixth `model' as a class.
\begin{itemize}
    \item \textbf{Example(s) 6: Congruent prisms arranged in columns.}  The triangular prism model, model 6a,  is depicted in Figure \ref{fig4}(b) and analysed in \cite{wc}. All cells are congruent triangular prisms
    with parallel longitudinal axes and front facets which are not aligned with each other. The prisms are arranged in columns. Similarly packed structures with all cells congruent and arranged in columns can be achieved with quadrilateral, pentagonal or hexagonal front facets (models 6b, 6c and 6d respectively). In each case, we confine our attention to models where the planar \tes, on a plane cut orthogonal to the column axis, is side-to-side and comprising, of course, congruent polygons for cells. All of these models have $\VE=4,\ \xi = \tf12$ and $\kappa = 0$, but the other parameters differ.
    \begin{itemize}
      \item[(a)] triangular prisms: $\EP = \tf92, \PV = \tf{27}{4}, \psi =5$ and $\tau =4$.
       \item[(b)] quadrilateral prisms: $\EP = \tf72, \PV = \tf{28}{5}, \psi =3$ and $\tau =2$.
       \item[(c)] pentagonal prisms (based on the so called Cairo side--to--side \tes\ of the plane with congruent pentagons, see Wells \cite{wel}, p.23): $\EP = \tf{16}{5}, \PV = \tf{16}{3}, \psi =\tf{12}{5}$ and $\tau =\tf75$.
       \item[(d)] hexagonal prisms: $\EP = 3, \PV = \tf{36}{7}, \psi =2$ and $\tau =1$.
    \end{itemize}
    These models belongs to a larger class of ``column \tes s" discussed in Section 5.
\end{itemize}
Then follow a number of examples constructed by adding dividing planes to some or all of the cells in Examples 1-6.
\begin{itemize}
     \item \textbf{Example 7:} A cube can be divided into three congruent
     pyramids. This is achieved with three congruent triangular plates
     which meet on one diagonal $AB$ of the cube, as shown in Figure \ref{fig4}(a). Of the six facets of the original cube, three remain as undivided squares (see those facets adjacent to $A$) whilst three are divided by a diagonal chord (those adjacent to $B$). The diagonal chosen for each cube could be randomised in some way, but here we choose a non-random way that produces a facet-to-facet \tes. Consider a larger cube containing eight of the original cubes. Now, starting from the centre of the large cube draw diagonals into each of the eight smaller cubes. Place three triangular plates in each cell meeting on this diagonal, choosing the method which puts the end $B$ (the corner where these triangles have their smallest angle) in the centre of the large cube. So, no parts of the $16$ added plates are visible from outside the large eight-celled cube because only the faces adjacent to $A$ are on the outside surface of this large cube. Then we tessellate $\R^3$ as a cubic lattice using the large cube. Here $(\VE, \EP)=(8,4)$ and $\PV=\tf{16}{5}$.
        \item \textbf{Example 8:} Each cell of a Delaunay \tes\ is a
        tetrahedron.  We independently divide each cell into two
        tetrahedra  by joining one of the cell's ridges, randomly chosen
        from the six it has, to a uniformly random point on the opposite
        ridge. The Delaunay is facet-to-facet, but the
        \emph{Divided Delaunay} is not. This model is analysed in full in \cite{wc} where it is shown that $\VE = \frac{14(5+16\pi^2)}{35+32\pi^2}\approx 6.50, \EP=\frac{72\pi^2(175+176\pi^2)}{7(5+16\pi^2)(35+24\pi^2)}\approx 4.38, \PV = \frac{9(175+176\pi^2)}{16(35+24\pi^2)}\approx 3.96, \xi = \frac{64\pi^2}{7(5+16\pi^2)}\approx 0.55, \kappa = 0, \psi = \frac{8\pi^2(35+528\pi^2))}{(35+24\pi^2)(35+32\pi^2)}\approx 4.34$ and $\tau = \frac{32\pi^2(102\pi^2 -35)}{(35+24\pi^2)(35+32\pi^2)}\approx 3.22$.
\end{itemize}

The construction in our next class of models starts with a planar \tes, so firstly a few remarks on notation will be apposite. The notations $\XY$ and $\la_X$ are also used in the planar case and only one adjacency, $\VE$, and one intensity, $\la_V$, are required to express all others. Only one `interior' parameter, namely $\phi$, is required (\cite{cow78},\cite{cow80}).
\begin{itemize}
    \item $\phi :=$ the proportion of vertices  which are contained in the interior of some cell-side. Note that $\phi = \mu_{VZ_1}\nc$.
\end{itemize}
These are called $\pi$-vertices because one angle between
consecutive edges emanating from the vertex is $\pi$. Planar \tes s
are \emph{side-to-side} if $\phi=0$, or not if $\phi>0$. Naturally
$0\leq \phi \leq 1$ and, as shown in \cite{wc}, $3\leq \VE \leq
6-2\phi$. ( This is a smaller constraint domain for $\VE$ than the
region $3\leq \VE \leq 6 $ given in \cite{km}; in that study, $\phi$
was not used.)
\begin{itemize}
    \item \textbf{Example(s) 9; the \emph{stratum} models:} Given any
    stationary \tes\ $\mathcal Y'$ in the plane, we can construct
    a \tes\ of $\R^3$ by firstly making each cell of $\mathcal Y'$ into a prism
    of depth $1$ unit by adding plates orthogonal to the plane in
    which $\mathcal Y'$ lies. This forms a stratum of depth $1$. We
    then place further copies of the stratum next to the first
    (without any offset of $\mathcal{Y'}$) and so on until $\R^3$ is filled. Mecke
    introduced this structure in \cite{mec84}, together with the
    name \emph{stratum}. Given a $(\VE', \phi)$-planar \tes, we can
    write the seven parameters ($\VE, \EP, ...)$ for the
    stratum model as functions of $\VE'$ and $\phi$.
\begin{alignat}{2}
    \VE =&\ \VE' +2 \qquad\qquad& \xi =&\ \frac{2\phi}{\VE'+2}\notag\\
    \EP =&\ \frac{6\VE'}{\VE'+2} & \psi =&\ 2\phi\qquad\qquad \tau = \phi \label{stratum}\\
    \PV =&\ \frac{3\VE'}{\VE'-1} & \kappa =&\ 0.\notag
\end{alignat}
Note that the three-dimensional \tes\ is \ftf\ if and only if the planar \tes\ is  side-to-side (which is equivalent to $\phi=0$).

The proof of these identities is deferred to a supportive paper
\cite{nwc}. We note here that all stratum models lie on the curve
$\EP = 6(1-2/\VE)$ and  with $5 \leq \VE \leq 8-\psi$. Dots
annotated $9a$--$9d$ on Figure \ref{fig1} correspondent to (a)-(d)
below.
\begin{itemize}
   \item[(a)] Stratum constructed from a Voronoi \tes\ in the plane, so $(\VE,\EP)=(5, \tf{18}{5})$ and $\PV=\tf92$.
   \item[(b)] Stratum constructed from the planar Delaunay \tes. Thus $(\VE,\EP)=(8, \tf{9}{2})$ and $\PV=\tf{18}{5}$.
   \item[(c)] Stratum constructed from the superposition of planar
   Voronoi \tes\ and its dual Delaunay \tes s (see \cite{wc}, Figure \ref{fig1}(a), where it is shown that $\VE'=4$ and $\phi=0$). So $(\VE,\EP)=(6, 4)$ and $\PV=4$.
   \item[(d)] Stratum constructed from the planar STIT \tes. Cyclic parameters are the same as for $9a$; also $\xi=\tf25, \psi=2, \tau=1$ and $\kappa=0$.
\end{itemize}
Dots $9a$ and $9b$ in Figure \ref{fig1} are the bounding cases.

   \item \textbf{Example(s) 10; the central-point models:} In this
   class of models, we commence with any stationary \tes\ of $\R^3$.
   In the interior of some cell, a point $Q$ is placed. Then, from
   every apex of that cell, a line segment is drawn to $Q$. Also, from
   every ridge of the cell a triangular plate is constructed with base
   being the ridge and $Q$ the opposite corner. Repeat this for all
   cells. Thus, every cell is itself partitioned into many pyramids, as
   many as there are facets of the original cell.

   The parameters are as follows.
  \begin{alignat*}{2}
  \VE  =&\ \frac{2 \PV' ( \VE'-4 - \VE'\EP'  + 2 \kappa' + 2 \psi')}{\PV' (\VE'-4)-\VE' \EP'  }\quad&\xi=&\ \frac{\VE' \xi'}{\VE'(\EP'-1)+4-2\kappa'-2\psi'}\\
   \EP=&\ \frac{4\VE'\EP'-3\VE'\xi'-4\psi'}{4 - \VE' + \VE'\EP' - 2 \kappa' - 2 \psi'}\quad&\kappa=&\ \frac{2\PV'\kappa'}{\VE' \EP'-\PV'(\VE'-4)}\\
       \PV=&\ \frac{\PV'(4\VE'\EP'-3\VE'\xi'-4\psi')}
    {\VE'\EP'(\PV'+1)-\PV'(\VE'\xi'+2\psi')}\quad&\psi=&\  \frac{4 \PV' \psi'}{\VE'\EP'-\PV'(\VE'-4)}\\
    &&\tau =&\ \frac{2 \PV' (\tau'+\psi')}{\VE'\EP'-\PV'(\VE'-4)}.
  \end{alignat*}

The dots marked $10a, 10b, ... $ in Figure \ref{fig1} are
constructed from the following models:
\begin{itemize}
   \item[(a)] a Voronoi \tes\ of $\R^3$ and yielding $(\VE,\EP)=(\tf{288\pi^2}{35+24\pi^2}, 4)$ and $\PV=\tf{576\pi^2}{7(5+24\pi^2)}$;
   \item[(b)] a  Delaunay \tes, yielding $(\VE,\EP)=(\tf{10(7+24\pi^2)}
   {35+24\pi^2}, \tf{576\pi^2}{7(5+24\pi^2)})$ and $\PV=3$;
   \item[(c)] the STIT \tes, with $(\VE,\EP)=(\tf{40}{7}, \tf{21}{5}), \PV=\tf{84}{19}, \xi= \tf35, \kappa = \tf47, \psi=\tf{24}{7}$ and $\tau = \tf{20}{7}$;
    \item[(d)] the cubic lattice \tes, with $(\VE,\EP)=(11, \tf{48}{11})$ and $\PV=\tf{16}{5}$;
    \item[(e)] the triangular prism \tes, defined in Example 6(a)  and yielding $(\VE,\EP)=(6, \tf{23}{4}), \PV=\tf{69}{13}, \xi= \tf14, \kappa = 0, \psi=\tf{15}{2}$ and $\tau = \tf{27}{4}$;
    \item[(f)]-- (i) further iterations of 6(a).
\end{itemize}
These iterations of Example 6(a) are connected by line segments in
Figure \ref{fig1}. Note that, because iteration decreases $\xi$, the
iterates have $(\VE, \EP)$ points which are moving toward the region
$\EP \leq 6(1-2/\VE)$ where \tes s with $\xi =0$ can exist. In the
limit, these iterations of $(\VE, \EP)$ reach $(8,\tf92)$ where
Example 9(b) is positioned.

Note that the constructed \tes s in Examples 10 are \ftf\ if and
only if the starting \tes\ is also \ftf.
\end{itemize}

\section*{3. Constraints in the \ftf\ case}

In this section we assume that the \tes\ is \ftf. Consider again the inequalities (\ref{simple}) and (\ref{simplex}). They apply to every cell of the \tes. So, $\ZV$ equals the average of $f_0$ over all cells. Likewise, $\ZE$ and $\ZP$ are the averages of $f_1$ and $f_2$ respectively. So (\ref{simple}) leads to
\begin{equation*}
    3 \ZV \leq 2\ZE,
\end{equation*}
which in turn (from the entries in Table 1) is equivalent to
\begin{align*}
    \frac{\PV f(2)}{f(\PV)} \leq&\ \frac{\VE\ \EP\ \PV}{f(\PV)},
\end{align*}
or to
\begin{equation}\label{fund}
    \EP \leq 6\Bigl(1-\frac{2}{\VE}\Bigr).
\end{equation}
This becomes a fundamental inequality for \ftf\ \tes s. The equality
in (\ref{fund}) holds if and only if all cells are \emph{simple}
polyhedra --- which means that each apex of the polyhedron must have
three ridges of the polyhedron emanating from it. The other
inequality, (\ref{simplex}), is less fruitful as it is equivalent to
something we know already, namely $\PV\geq 3$.

\begin{figure}[ht]
    \psfrag{EP}{{ \scriptsize$\EP$}}
    \psfrag{PV}{{ \scriptsize$\hspace{-2mm}\PV$}}
    \begin{center}
    \includegraphics[width=165mm]{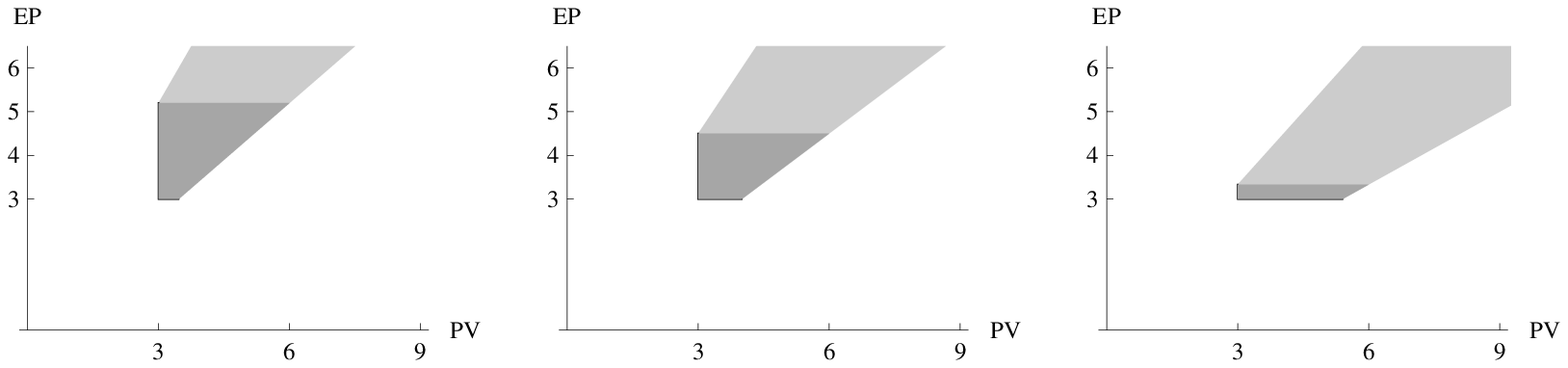}
      \\
          {     \scriptsize (a) $\VE =15$\hspace{3.8cm}(b) $\VE=8$ \hspace{3.8cm}(c) $\VE=\tf92$}
        \caption{\label{fig5}\scriptsize The darkly shaded region is the range of $(\PV,\EP)$ for various
        $\VE$ when the \tes\ is \ftf. Note that the upper and lower
        boundaries of the dark zone coincide when $\VE=4$.  We shall see
        later, in Section 6, the theory behind the light grey
        region. This region,
       an open set bounded on the left by the sloping line
        $\EP=2(1-2/\VE)\PV$, is the range for non \ftf\ \tes s which
        have $\tau>0$. The sloping line which bounds both
        shaded regions on the right is
        $\EP=(1-2/\VE) \PV$. Tessellations which are not \ftf\ and have
        $\tau =0$
        are also in the darkly-shaded region (see Remark 9).
       The non \ftf\ theory is described later in
        the paper.}
    \end{center}
\end{figure}

Turning now to the rest of Table 1, we can systematically write
inequalities attached to many of the entries --- for example, $\ZV
\geq 4, \ZE\geq 6,\ \ZP\geq 4$ and $f(\PV)> 0$. The first three of
these provide lower bounds for $\PV$, but all of these are weaker
than $\PV\geq 3$ when (\ref{fund}) holds --- as it does when the
\tes\ is \ftf.  However, $\la_Z>0 \Llra  f(\PV)>0$ which is
equivalent to
\begin{equation}\label{PVupper}
    \PV <  \frac{\VE\EP}{\VE-2},
\end{equation}
this upper bound for $\PV$ being $\le 6$ in the \ftf\ case because
of (\ref{fund}).  We have thus completed the proof of
(\ref{case1th}). The darkly shaded regions of Figure \ref{fig5} show
the permitted range for $(\PV,\EP)$ for various $\VE$, illustrating
the inequalities in (\ref{case1th}).

\section*{4. Facet-to-facet examples with large $\VE$ or large $\ZV$}

Figure \ref{fig6} shows how some parameters in a stationary \tes\ of
$\R^3$ can be arbitrarily large. This is in contrast with the
situation for planar \tes s. Both \ref{fig6}(a) and \ref{fig6}(b)
are based on the cubic lattice \tes, with structure added in the
interior of each cube.
\begin{itemize}
\item \textbf{Example 11; model with unbounded $\VE$:} In Figure \ref{fig6}(a), a
vertical line joins the centres of the top and bottom facets of the
cube. Called these centres $T$ and $B$. This line $TB$ is further
divided into $(k+1)$ edges by the addition of $k$ equally-spaced
vertices $V_1, V_2, ..., V_k$. The case $k=1$ is shown in Figure
\ref{fig6}(a). Around the boundary of the bottom facet are placed
$4n$ vertices in addition to those at the face's corners, $n$ in the
interior of each side of the square facet. Figure \ref{fig6}(a)
shows the case $n=2$. Similarly, the top facet has these added
vertices on its boundary.
\begin{figure}[ht]
    \psfrag{T}{$T$}
    \psfrag{B}{$B$}
    \psfrag{V}{$V_1$}
    \begin{center}
    \includegraphics[width=68mm]{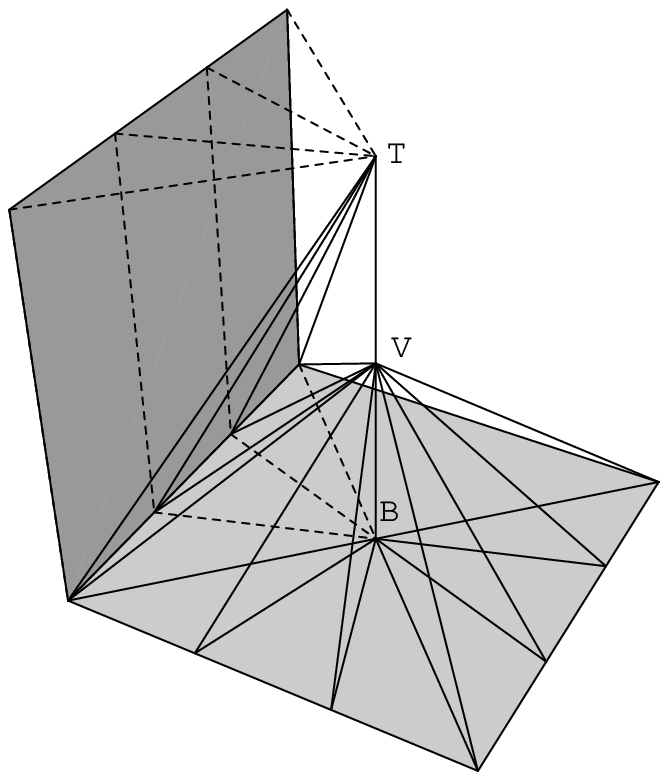}\hspace{0cm}
    \includegraphics[width=65mm]{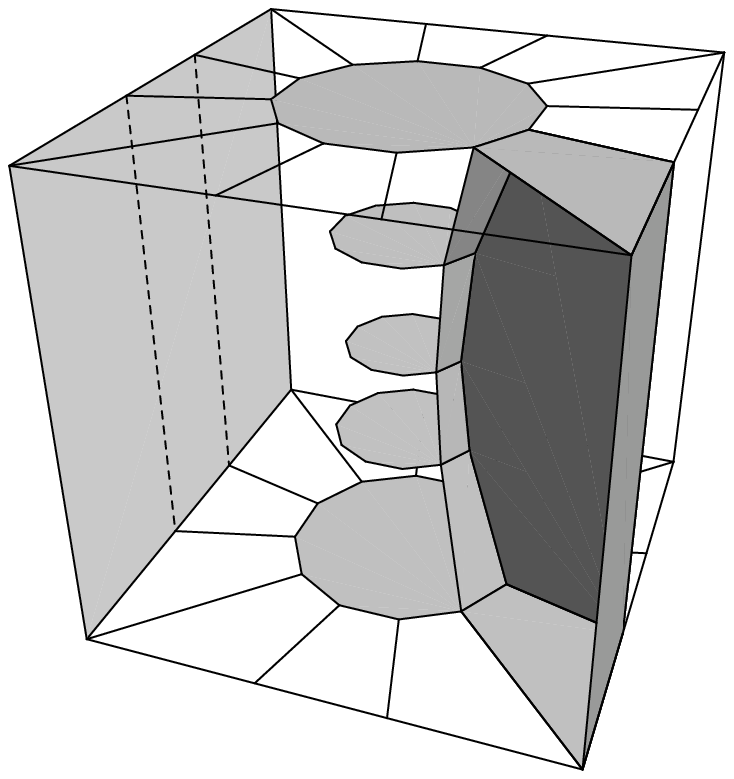} \\
                 \renewcommand{\baselinestretch}{.1}\small
          {\scriptsize (a)\hspace{8cm}(b)}
        \caption{\label{fig6} \scriptsize (a) Model with unbounded $\VE$; (b)
        model with unbounded $\ZV$.}
    \end{center}
\end{figure}
From each vertex $V_j, j=1,2,...,k$, and $T$ and $B$ too, edges are
drawn to the $4(n+1)$ vertices on the boundary of the bottom facet
(as shown). Edges also join $T$ with the $4(n+1)$ vertices on the
boundary of the top facet. Finally $T$ and $B$ are joined, creating
many edges because the joining line-segment passes though $V_j,
j=1,2,...,k$. Plates are now formed in the obvious places seen in
the drawing. A facet-to-facet \tes\ results.

It can be shown that:
\begin{equation*}
    \VE = \frac{2 (12 + 5 k + 12 n + 4 k n)}{2 + k + 2 n};\quad
    \VP = \frac{8 (7 + 3 k) (1 + n)}{2 + k + 2 n};\quad
    \VZ = \frac{4 (9 + 4 k) (1 + n)}{2 + k + 2 n}.
\end{equation*}
So, $\VE$ can be made arbitrarily large (and this reflects in the
other mean adjacencies $\VP$ and $\VZ$ associated with the typical
vertex). All other adjacencies are bounded as $k$ and/or $n$ becomes
large. For example:
\begin{equation*}
    \EP= \frac{8 (7 + 3 k) (1 + n)}{12 + 5 k + 12 n + 4 k n}; \qquad\qquad\qquad
    \PV = \frac{4 (7 + 3 k)}{9 + 4 k}.
\end{equation*}
Some points for small $n$ and $k$ appear on the plot of Figure
\ref{fig1}, annotated $(k,n)=(2,0)$ as `11a', $(k,n)=(0,0)$ as `11b'
and $(k,n)=(2,1)$ as `11c'. For larger $k$ and $n$, the points lie
off the graph to the right, but always below the curve (approaching
the curve as $k\rightarrow \i$).

\item \textbf{Example 12; model with unbounded $\ZV$:} In Figure \ref{fig6}(b),
there are also $4(n+1)$ vertices on the boundary of both the bottom
and top facets of the cube. There is a central core made up of $k+1$
prisms, each bounded above and below by a $4(n+1)$-sided polygonal
plate. Other cells are constructed outside the central core as per
the illustration (which uses $n=2$ and $k=3$). The \tes\ is \ftf,
with:
\begin{equation*}
    \ZV = \frac{8 (5 + 2 k) (1 + n)}{5 + k + 4 n};\quad
    \ZE = \frac{12 (5 + 2 k) (1 + n)}{5 + k + 4 n};\quad
    \ZP = \frac{2 (15 + 5 k + 14 n + 4 k n)}{5 + k + 4 n}.
\end{equation*}
We note that these adjacencies can be made arbitrarily large. Our
cyclic trio are, however, bounded above.
\begin{align*}
        \VE =&\ \frac{2 (15 + 8 k + 16 n + 8 k n)}{5 + 4 k + 6 n + 4 k n};\\
    \EP=&\ \frac{12 (5 + 2 k) (1 + n)}{15 + 8 k + 16 n + 8 k n}; \\
    \PV =&\ \frac{12 (5 + 2 k) (1 + n)}{15 + 5 k + 14 n + 4 k n}.
\end{align*}
These formulae lead to a datum for Figure \ref{fig1}, for the case
$n=k=3$. All such data for this example would lie on the fundamental
curve, with $4 \leq \VE \leq 6$.
\end{itemize}
The examples above in this section show rather extreme behaviour,
demonstrating the unbounded character of $\VE$ and $\ZV$ ---
something already noted for \ftf\ tilings by Ziegler at a 2002
conference (reported in \cite{zie}).

Example 11 is especially useful to our Figure \ref{fig1}, as it
shows that valid \tes s exist far to the right of the diagram, below
the curve. It is also useful in a rather esoteric way, in
\emph{mixtures of models}.

\section*{5. Tessellations which are mixtures}

Formally, a random \tes\ is
    a mapping from a probability space $(\Omega, \mathcal{T}, \Pr)$
    into the space of all \tes s of the holding space, here $\R^3$.
    The probability measure $\Pr$ can be of the form
    $\al \Pr_1+ (1-\al)\Pr_2$ where $0\leq \al\leq 1$, $\Pr_1$ and
    $\Pr_2$ being
    probability measures. In practical terms, this means that
    with probability $\al$ one \tes\ model is realised and with
    probability $(1-\al)$ it is another model. As $\al$ varies through
    the interval $[0,1]$, a curve is traced between the points
    of Figure \ref{fig1} which represent $\Pr_1$ and $\Pr_2$ and all mixtures
    of these lie on the curve.

\begin{itemize}
    \item \textbf{Example(s) 13; Two mixture \tes s and some mixture curves.}  For example, the dashed curves in
    Figure \ref{fig1} show the mixture curves of $8$ and $7$,
    $10e$ and $2$ and $6a$ with model $11$ (having $k=n=100$ and
    therefore its $(\VE, \EP)$ point located far to the
    right of the region shown and a very small distance below
    the solid curve).
    Both the dots marked $13a$ and $13b$  are mixtures of
    $2$ (as $\Pr_1$) and $10e$ (as $\Pr_2$) with
    $\al = 0.8$. In $13b$ the vertex intensities of the two pure
    models are equal, whilst in $13a$ the vertex intensity of $10e$
    is twice that of $2$.
\end{itemize}
\textbf{Theory of mixture curves:} When we mix two \tes s of $\R^3$,
we can calculate the properties of the mixture from the following
two identities (using suitable choices for the object classes $X$
and $Y$). Using a superscript $[j]$ (note the square brackets) to
indicate that $\Pr_j$ is operative, we can write the following:
\begin{align*}
    \la_X =&\ \al \la_X^{[1]} +(1-\al)\la_X^{[2]};\\
    \XY=&\ \frac{\la_X^{[1]} \al \XY^{[1]} +\la_X^{[2]} (1-\al)
    \XY^{[2]}}{\al \la_X^{[1]} +(1-\al)\la_X^{[2]}}.
\end{align*}
The data points in Figure \ref{fig1} are $(\VE, \EP)$. If such a
point is a mixture, we can write:
\begin{align*}
    \VE=&\ \frac{\la_V^{[1]} \al \VE^{[1]} +\la_V^{[2]} (1-\al)
    \VE^{[2]}}{\al \la_V^{[1]} +(1-\al)\la_V^{[2]}}; \\
    \EP=&\ \frac{\la_E^{[1]} \al \EP^{[1]} +\la_E^{[2]} (1-\al)
    \EP^{[2]}}{\al \la_E^{[1]} +(1-\al)\la_E^{[2]}}.
\end{align*}
As $\al$ moves through the interval $[0,1]$, a continuous curve is
traced from $(\VE^{[2]},\EP^{[2]})$ to $(\VE^{[1]},\EP^{[1]})$.
Eliminating $\al$ gives the curve ($\EP$ versus $\VE$) as
\begin{equation}
\EP = \frac{\EP^{[1]} \VE^{[1]} - \EP^{[2]}
\VE^{[2]}}{\VE^{[1]}-\VE^{[2]}} - \frac{(\EP^{[1]} - \EP^{[2]})
\VE^{[1]}\VE^{[2]}}{\VE(\VE^{[1]}-\VE^{[2]})}\label{join}
\end{equation}
when $\VE^{[1]}\ne \VE^{[2]}$. When $\VE^{[1]} = \VE^{[2]}$, the
curve is a vertical line-segment joining the two points.

Note that the result in (\ref{join}) is valid whether $\la_V^{[1]} =
\la_V^{[2]}$ or not. Also note that, if both $(\VE^{[1]},\EP^{[1]})$
and $(\VE^{[2]},\EP^{[2]})$ lie on a curve of the form $\EP =
A-B/\VE$, where $A$ and $B$ are constants, (\ref{join}) lies on the
same curve as shown by the following calculation.
\begin{align*}
    \EP =&\ \frac{(A-B/\VE^{[1]}) \VE^{[1]} - (A-B/\VE^{[2]}) \VE^{[2]}}{\VE^{[1]}-\VE^{[2]}}
    - \frac{[(A-B/\VE^{[1]}) -
(A-B/\VE^{[2]})]\ \VE^{[1]}\VE^{[2]}}{\VE(\VE^{[1]}-\VE^{[2]})}\\
    =&\ \frac{A\VE^{[1]} - A\VE^{[2]} -B+B}{\VE^{[1]}-\VE^{[2]}}
    - \frac{
B\VE^{[1]}-B\VE^{[2]}}{\VE(\VE^{[1]}-\VE^{[2]})}\\
    =&\ A \bigl(1- \frac{B}{\VE}\bigr).
\end{align*}
Our fundamental curve is of this form, with $A =6$ and $B=12$, and is therefore closed under the mixture operation.

{\sc Remark 3:} \emph{Readers accustomed to non-random tilings may
regard such mixtures of models as fraudulent. To them, the more
interesting question is ``what is the existence domain in Figure
\ref{fig1} using only non-mixture models?''. We do not yet have an
answer to this difficult question.}

\section*{6. Introductory theory for non \ftf\ \tes s}

\textbf{Proof of (\ref{case2}) for non \ftf\ \tes s:} The most
obvious change from the \ftf\ theory that we dealt with in Section 3
is that (\ref{fund}), the upper bound on $\EP$,  is no longer valid.
We can however retain (\ref{PVupper}), the upper bound on $\PV$,
because the argument for it remains sound in the non \ftf\ case. It
cannot now be argued, however, that this upper bound on $\PV$ is
$\leq 6$; the former argument had used (\ref{fund}).

We now address the lower bound for $\PV$ in the non \ftf\ case. To
find this bound, we initially consider the entity $\ZV$, the
expected number of \tes\ vertices adjacent to a typical cell. For
$\ZV$ to equal $4$, it is necessary that all cells be tetrahedra
--- and for there to be no \pies, because the vertices at the
termini of any \pies\ would raise $\ZV$ above $4$. So $\ZV =4 \Lra
\xi=0$. Therefore a non \ftf\ \tes\ has $\ZV>4$; we use this fact to
establish the bound (\ref{ZV4}) that follows. The entry for $\ZV$ in
Table 1 and the definition of $f$ in (\ref{fx}) are also used.
\begin{align}
    \mathrm{non\ facet}\hspace{-.7mm}-\hspace{-.7mm}\mathrm{to}
    \hspace{-.7mm}-\hspace{-.7mm}\mathrm{facet} \Lra  \ZV> 4
    \Llra& \ds{\frac{\PV f(2)}{f(\PV)}}> 4 \notag\\
    \Llra& \ds{\frac{\PV (\VE\EP - 2(\VE-2))}{\VE\ \EP-\PV(\VE-2)}}> 4 \notag\\
    \Llra& \frac{\VE\EP}{2(\VE-2)} < \PV.\label{ZV4}
\end{align}
This new lower bound for $\PV$ is $> 3$ if and only if
$\EP>6(1-2/\VE)$, so it plays no role in \ftf\ \tes s that have
$\ZV>4$. The inequality $\PV \geq 3$ remains operative for all \tes
s, \ftf\ or not, when $\EP<6(1-2/\VE)$. Because
\begin{align*}
    \EP =&\ 6\bigl(1-\frac{2}{\VE}\bigr) \Llra \frac{\VE\EP}{2(\VE-2)} = 3,
\end{align*}
(\ref{ZV4}) implies that $\PV > 3$ in the non \ftf\ case when
$\EP=6(1-2/\VE)$ (though the inequality remains $\PV\geq3$ in the
\ftf\ case). Thus (\ref{case2}) in the main theorem is proven. An
illustration of these inequalities is given in Figure \ref{fig5}.

{\sc Remark 4:} \emph{When $\EP>6(1-2/\VE)$, the lower bound $\tf12\VE\EP/(\VE-2)$ is also greater than lower bounds that can be similarly derived from $\ZE > 6$ and $\ZP> 4$.}

 \textbf{A focus  on $\EP$:} We have seen above that
our earlier proof of (\ref{fund}), which says that $\EP\leq
6(1-2/\VE)$ in the \ftf\ case, no longer holds. Figure \ref{fig1},
with its many black dots above the fundamental curve, goes further
and demonstrates that $\EP$ can be greater than $6(1-2/\VE)$ in the
non \ftf\ case. The following class of models reveal that $\EP$ is
indeed unbounded when $\VE=4$.

\begin{itemize}
    \item \textbf{Example(s) 14: Column \tes s.} Given a stationary planar tessellation ${\mathcal Y}'$ we construct a spatial tessellation as follows. Each cell of ${\mathcal Y}'$ is the base of an infinite cylinder perpendicular to the plane ${\mathcal E}$ in which ${\mathcal Y}'$ lies. These cylinders are called columns. They pack to fill ${\R}^3$. Now any column is intersected by planes parallel to ${\mathcal E}$ with constant separation $1$. The position of these cuts is stationary and independent of those in the neighboring columns (unlike the positioning in the \emph{stratum} models which aligned the cuts in all columns).

        The result is a stationary spatial tessellation --- the \emph{column} tessellation ${\mathcal Y}$. Any cell of ${\mathcal Y}$ is a right prism with height $1$ and a base facet which is a translate of a cell of ${\mathcal Y}'$. Due to the independence of cuts in the columns, no column tessellations are facet-to-facet.

        The topology of ${\mathcal Y}$ is determined by the topology of
        ${\mathcal Y}'$. Hence the seven parameters for a column tessellation
        are functions of some topological parameters of the  planar
        tessellation. Besides $\mu'_{VE}$ and $\phi$ we need
\begin{tabbing}
$\hspace{.5cm}\circ\ \mu'_{E\,V[\pi]} :=$ the mean number of $\pi$-vertices adjacent to the typical edge and\\[2mm]
$\hspace{.5cm}\circ\ \mu'^{(2)}_{VE}:=$ the second moment of the
number of edges adjacent to the typical vertex.
\end{tabbing}
We obtain
$$\begin{array}{rclrcl}
\mu_{VE} &=& 4 \label{VE} &
\xi &=& \frac{1}{2} + \frac{1}{4} \mu'_{E\,V[\pi]} \label{xi}\\[2mm]
\mu_{EP} &=& \ds{\frac{1}{2 \mu'_{VE}}}\,( 3 \mu'_{VE} + \mu'^{(2)}_{VE}) \label{EP}
 \qquad &
\kappa &=& \frac{1}{2}  \mu'_{E\,V[\pi]} - \ds{\frac{\phi}{\mu'_{VE}}} \label{kappa}\\[3mm]
\mu_{PV} &=& \ds{\frac{2}{3\mu'_{VE}-2}}\, ( 3 \mu'_{VE} + \mu'^{(2)}_{VE}) \label{PV}\qquad
 &
\psi &=& \ds{\frac{\mu'^{(2)}_{VE} + 3 \phi}{\mu'_{VE}}} -1-\frac{1}{2}  \mu'_{E\,V[\pi]} \label{psi} \\[3mm]
&&&\tau &=& \ds{\frac{\mu'^{(2)}_{VE} + \phi}{\mu'_{VE}}}- 2 \label{tau1}.
\end{array}$$
The proof of these identities is available in \cite{nwc}. Their use
can be seen by looking again at Examples 6(a)--6(d). As an exercise,
one can show that Example 6(c) obeys these identities; by using the
two-dimensional entities $\VE' = \tf{10}{3},
\mu'^{(2)}_{VE}=\tf{34}{3}$ and $\phi = \mu'_{E\,V[\pi]} = 0$, one
recovers the cited three-dimensional answers.

The following column \tes\ is based on a planar tessellation whose
second moment $\mu'^{(2)}_{VE}$  is  unbounded and, from that
observation, we establish that the spatial model has no upper bound
for $\EP$. The starting point is a stationary tessellation in
${\R}^2$, where all cells are unit squares which are positioned
side-to-side. Now each cell is divided by a further vertex in its
relative interior, the vertex having $4n$ emanating edges, $n$ of
them to each side of the square, $n \geq 1$. Moreover those interior
edges of two squares having a common side are disjoint; they do not
meet at a common vertex in the interior of that common side. Hence
any square side has $2n$ new vertices in its relative interior ---
and they are all $\pi$-vertices with three emanating edges.

For such a planar tessellation we obtain
$$\begin{array}{cclrcl}
\mu'_{VE} &=& \ds{\frac{1\cdot 4 + 2 \cdot 2n \cdot 3 + 1 \cdot 4n}{1 + 2
\cdot 2n + 1}} \  = \ \ds{\frac{2(4n+1)}{2n+1}}&\qquad\quad
 \phi &=& \ds{\frac{2n}{2n+1}} \\[3mm]
\mu'^{(2)}_{VE} &=&  \ds{\frac{1\cdot 16 + 2 \cdot 2n \cdot 9 + 1
\cdot 16n^2}{1 + 2 \cdot 2n + 1}} \  = \
\ds{\frac{2(4n^2+9n+4)}{2n+1}}&\qquad\quad \mu'_{E\,V[\pi]} &=&
\ds{\frac{6n}{4n+1}}.
\end{array} $$

Hence $\mu_{EP}$ for column tessellations is unbounded and all the
tessellations from that model lie on the line $\mu_{VE}=4$ with
$\mu_{EP} \geq 4$ (see Figure \ref{fig1} where the point which
represents Example 14 has $n=8$ and $(\VE, \EP) =
(4,\tf{431}{66})$).

\end{itemize}

{\sc Remark 5}: \emph{Note that every point $Q$ on the curve
$\EP=6(1-2/\VE)$ represents some \tes. This is obvious from our
earlier theory of mixture curves if $Q$ lies on the curve between
 the three \tes s represented by $(4,3)$ and some
point $R$ on the curve. For illustration, taking $R$ as the point
from Example 2, we conclude that every point $Q$ on the curve
between $(4,3)$ and $R$ represents a \tes\ because the curve is of
the generic form $\EP = A - B/\VE$ closed under the mixture
operation. There is always such a $R$, however, for every such $Q$.
This is so, because for any $R$ higher on the curve than $(\VE,\EP)$
of Example 2, can be achieved by mixture. One can always find a
column \tes\ represented by $(4,\EP^{[1]})$ for some suitably large
$\EP^{[1]}$ and a \tes\ based on Example 11 with a suitably large
$\VE^{[2]}$ such that their mixture curve passes through $R$.}

{\sc Remark 6}: \emph{Since $\EP$ is unbounded when $\VE=4$, we can
use the mixture concept and Remark 5 to establish that $\EP$ is
unbounded for all $\VE$.}

\textbf{When $\xi>0$, other interior parameters can be zero:} For
 \ftf\ \tes s, we know that $\xi = \kappa =
\psi=\tau=0$. We also know that $\xi>0$ if and only if the \tes\ is
\emph{not} \ftf. The other interior parameters, $\kappa, \psi$ and
$\tau$, can still be zero when $\xi>0$, as the following two
examples show.

\begin{itemize}
    \item \textbf{Example 15:}
    Start with the cubic lattice of Example 5, aligned to the three
    Cartesian axes.  Divide each
    cube into two congruent triangular prisms using a rectangular
    plate. Vary the orientation of the plate, as shown in Figure
    \ref{fig4}(c). Clearly every such plate creates two $\pi$-edges in the
    resulting \tes. It is easy to show that $(\VE,\EP)=(10,4),
    \ \PV=\tfrac{10}{3},\ \xi=\tfrac25$ and $\kappa=\psi=\tau=0$.
    \item \textbf{Example 16:}
    Also start with the aligned cubic lattice.
    Partition each cube into three congruent pyramids as per
    Figure \ref{fig4}(a), with all of the $AB$ chords (one for each cube)
    being parallel and with all points marked $A$ being at the low end
    of the chord.
    Thus each cube-facet with a dividing diagonal is next to
    one with no such diagonal. These diagonals are therefore $\pi$-edges
    in the new \tes. As they represent $\tf37$ of all edges, $\xi=\tf37$.
    No vertices lie in the interior of any facet, ridge or plate--side,
    so $\kappa = \psi=\tau=0$. For the record, $(\VE,\EP)=(14,\tf{27}{7})$
    and $\PV=3$.
\end{itemize}
When $\xi>0$ it is possible that $\psi>0$ and $\kappa=0$
(Examples 6, 8 and 9) or that $\kappa>0$ and $\psi=0$ (Example 17
below).
\begin{itemize}
    \item \textbf{Example 17:}
    Adjust each column in Example 6a so that
    the triangular facets are aligned on parallel
    planes. There is now no offset between columns and we have a stratum
    \tes.
    In every second stratum, partition every
    cell into three congruent prisms as shown in Figure \ref{fig7}(a).
    The new vertices created are all hemi-vertices and one can show
    that $\kappa = \tfrac23$ and $\psi=\tau=0$. Also $\xi=\tfrac{6}{11},
    \ (\VE,\EP) =(\tfrac{22}{3},\frac{42}{11})$ and $\PV=\frac72$.
\end{itemize}

\section*{7. Books with spines}

Three of the interior parameters, $\kappa, \psi$ and $\tau$, are
related to vertices being in the interior of ridges, plate--sides or
facets, whereas the other interior parameter $\xi$ is related to
\pies. To detect constraints for the interior parameters it makes
sense to utilize relations between $\pi$-edges and interior vertices
based on the combinatorial topology of the tessellation. There is a
wide variety of arrangements and, to handle all those which are of
interest, we use a new concept: \emph{books} with their
\emph{spines}. We introduce the concept by explaining how
plate--sides whose interior contains some vertex $v$ can form the
\emph{pages} of a \emph{book}.

{\sc Definition 1:} \emph{Let $v$ be a vertex and $\ell$ be a
line passing through $v$. So $v\subset \ell$. We define a
\emph{page} as a plate which has one of its sides $\subset \ell$
 and the interior of that side $\supset v$. A \emph{book-cover} is a
 cell-facet which has  one of its sides $\supset v$ and contains
$k\geq 2$ plates that have a corner coinciding
 with $v$. Note: Exactly two of these plates have a side that is
 $\subset \ell$. Furthermore, the facet contains
 $(k-1)$ $\pi$-edges that emanate from $v$  and it may
 have other plates and \pies\ that are not adjacent to $v$. }

\vspace{-.3cm}\begin{figure}[ht] \psfrag{P1}{\small
\hspace{-2mm}$P_1$} \psfrag{C1}{\small\hspace{-2mm}$C_1$}
\psfrag{C2}{\small \hspace{-2mm}$C_2$}\psfrag{P2}{\small
$P_2$}\psfrag{P3}{\small $P_3$}
\begin{center}
    \includegraphics[width=40mm]{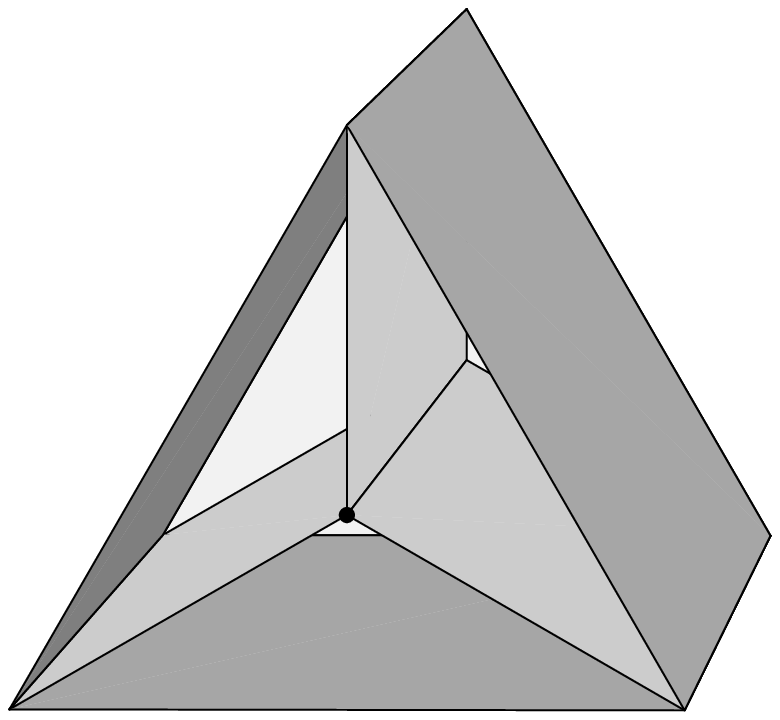}\hspace{-1.0cm}
    \includegraphics[width=77mm]{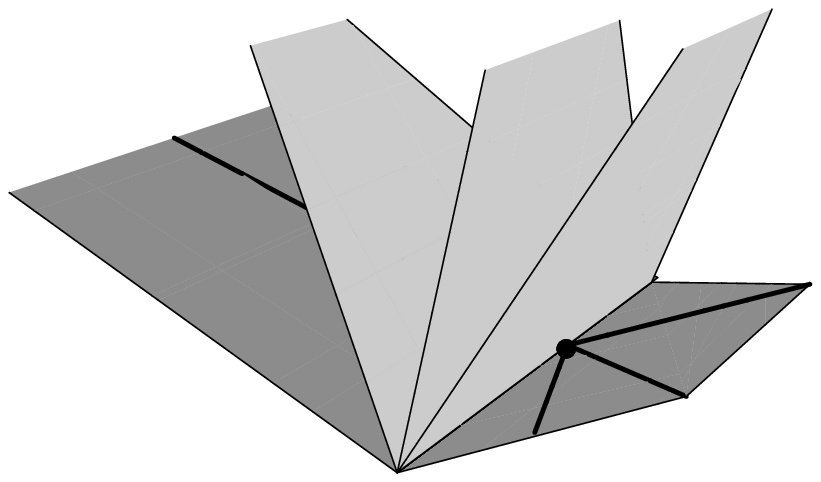}\hspace{-2.1cm}
   \includegraphics[width=80mm]{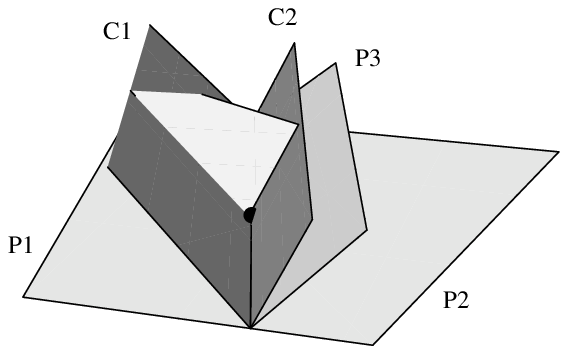}\\\vspace{-.8cm}
          {\scriptsize (a) \hspace{4cm}
          (b)\ \  \includegraphics[width=15mm]{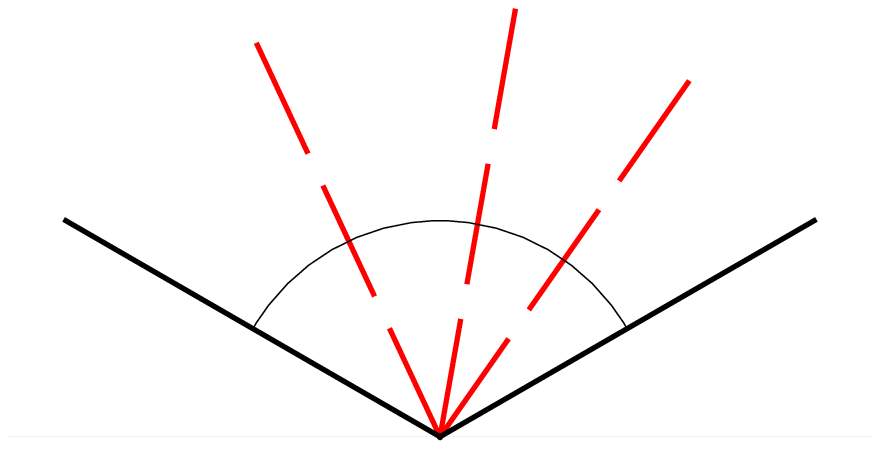}\hspace{4cm}
          (c)\ \ \includegraphics[width=15mm]{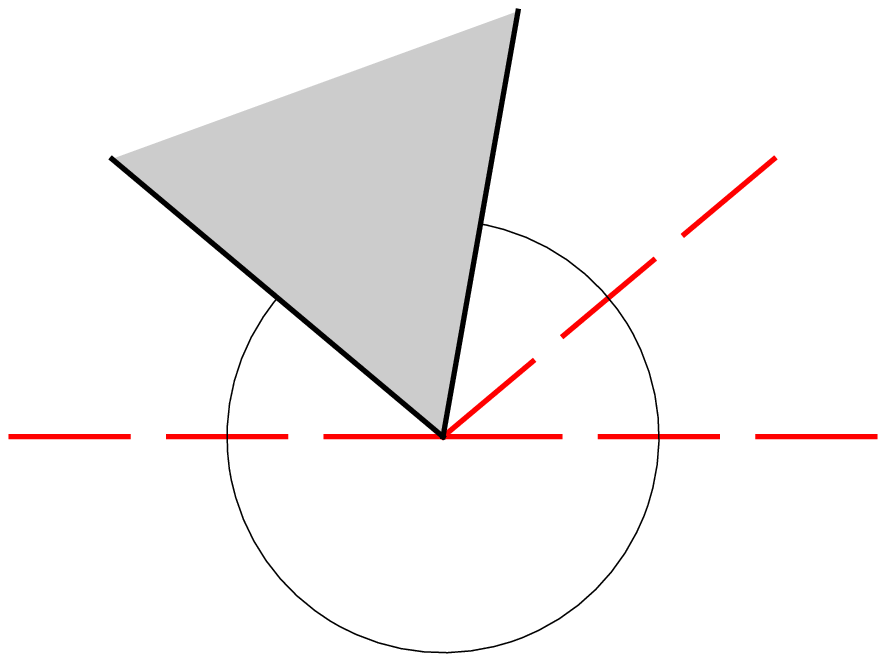}   }
        \caption{\scriptsize  \label{fig7} (a) All cells in every
        second stratum of Example 17  are partitioned as shown,
        into three congruent triangular prisms. (b) A $3$-page book
        with its two book-covers in darker shading. There are always
        edges on the book-covers emanating from the vertex which
        is on the book's `spine' (c)  A vertex adjacent to three plate--side-interiors,
        three ridge-interiors and one facet-interior (the horizontal plates
        $P_1$ and $P_2$ forming in union a facet of a cell below).
        Note: if the light-coloured plate between the two dark structures $C_1$ and $C_2$ and with a corner touching the vertex
        were removed,
        then the vertex would not exist. Using the `book' terminology, we see \emph{pages} $P_1, P_2$ and $P_3$ and two \emph{bookcovers} $C_1$ and $C_2$.
        }
\end{center}
\end{figure}

\vspace{0cm} \noindent In Figure \ref{fig7}(b), three pages  and two
darkly-shaded book-covers (the fully visible one having $k=4$) are
shown, in a neighbourhood of a vertex. The overall appearance is
like a book.

{\sc Definition 2:} \emph{At vertex $v$, a \emph{book} with $p\geq 0$
pages is a collection $\mathcal{B}$ comprising $p$ pages and two
book-covers in such a way that a circular arc, with small radius and
centred at $v$, can be drawn from one bookcover to the other
encountering no plates other than those which form the $p$ pages.
The line $\ell$ mentioned above is called the \emph{spine}
of the book.}

\vspace{0cm}
\begin{figure}[ht] \psfrag{P1}{\small\hspace{-2mm}
$P_1$} \psfrag{C1}{\small\hspace{-1mm}$C_1$} \psfrag{C2}{\small
\hspace{-2mm}$C_2$} \psfrag{P2}{\small $P_2$}
\psfrag{C3}{\small$C_3$}
    \begin{center}
    \includegraphics[width=52mm]{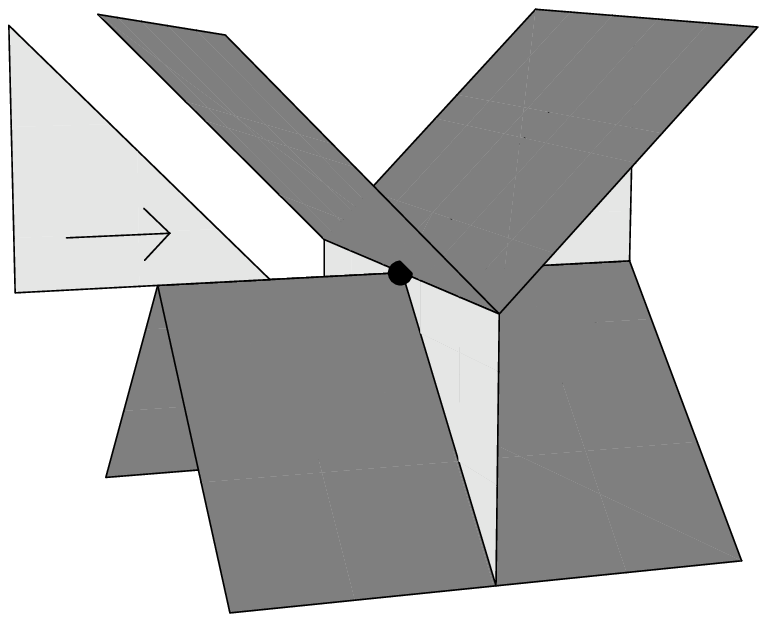}\hspace{-1.8cm}
    \includegraphics[width=88mm]{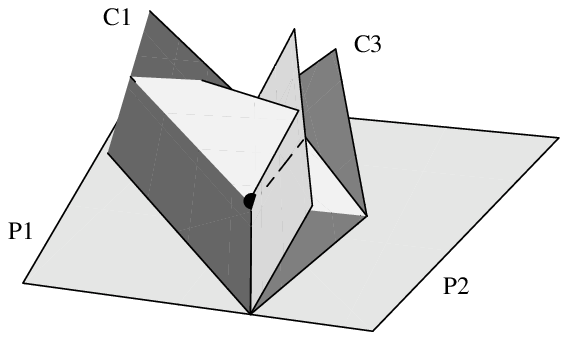}\hspace{.2cm}
     \includegraphics[width=40mm]{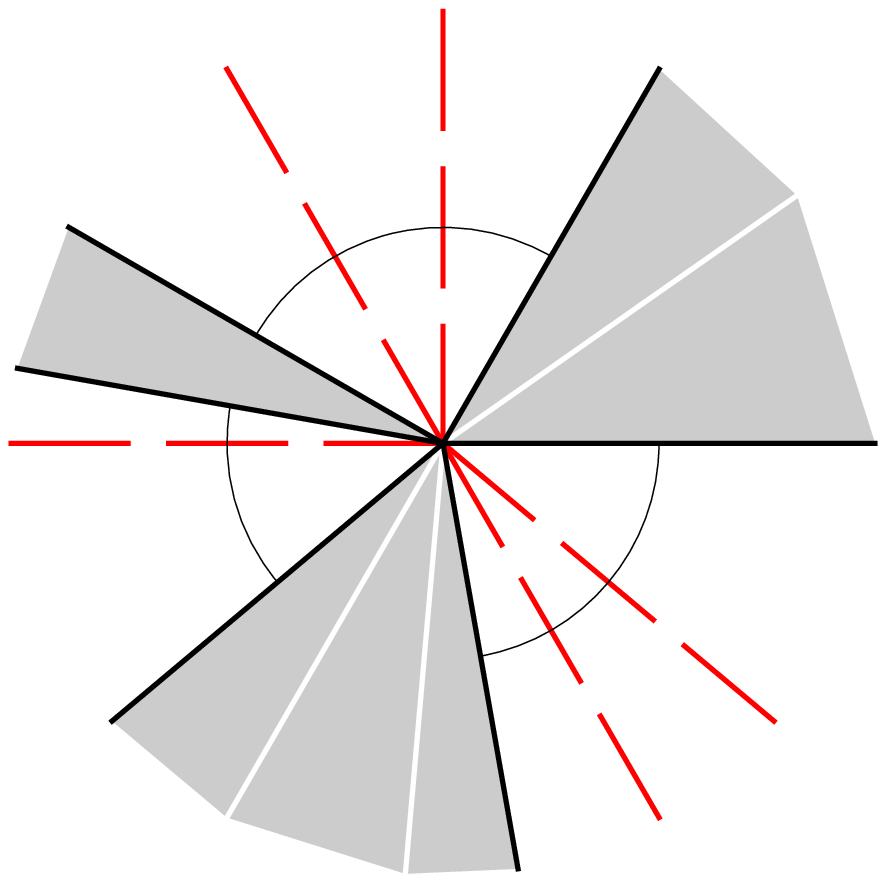}\\\vspace{-.5cm}
     {\scriptsize (a) \ \  \includegraphics[width=15mm]{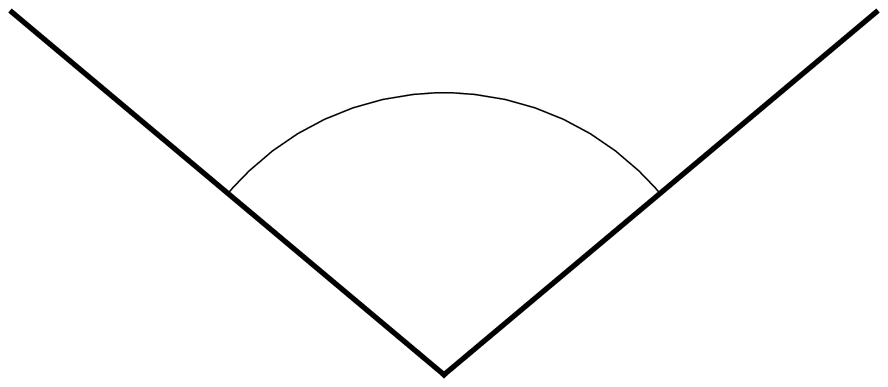}
      and \ \  \includegraphics[width=15mm]{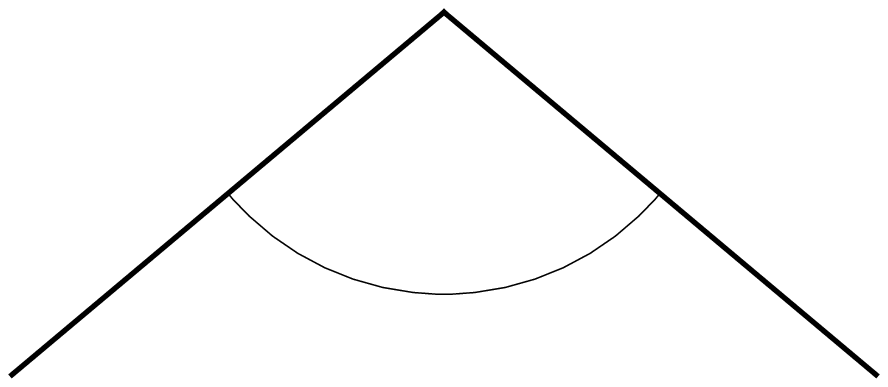}\hspace{2cm}(b) \ \  \includegraphics[width=15mm]{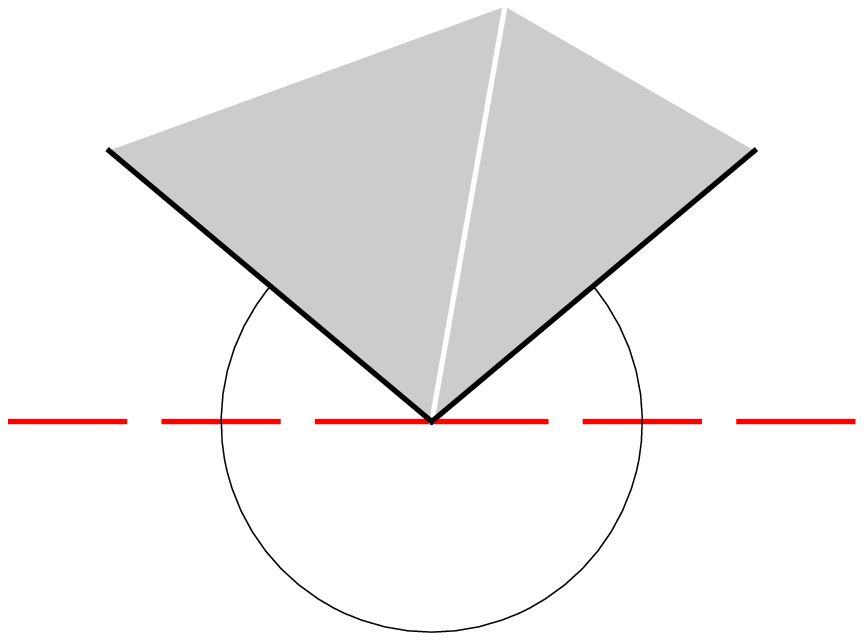} \hspace{5cm}(c)}
        \caption{\scriptsize  \label{fig8} (a) Two-books associated with a common
        vertex and with different
        spines. (b) A hemi-vertex $v$ in the interior of a  facet made from the coplanar pages $P_1$ and $P_2$. There is one book, with two pages. The additional architecture  has two plates with corners touching $v$ and both of these plates have created a \pie\ on the grey structure that follows the spine. (c) The schematic diagram for a vertex having three books each with a common spine. }
    \end{center}
\end{figure}

\vspace{0cm} {\sc Remark 7:} \emph{It is time to describe the
schematic figures, which look at a vertex $v$ along the spine. The
bookcovers are black whilst the pages are  dashed and red. The
circular arc is shown and shaded sectors represent the plate (or
plates) that have a corner on $v$, thus ensuring that $v$ is indeed
a vertex. Figure \ref{fig7}(c) provides a simple example. In Figure
\ref{fig8}(c), which schematically shows three books with a common
spine, we see shaded sectors sub-divided by  white lines. Such a
white line, also evident in Figure \ref{fig8}(b), represents a
coplanar structure that `lies along' the spine in a neighbourhood of
$v$, but a structure that is neither a page nor a bookcover.}

If $v$ is a hemi-vertex, then there may be a book with two
\emph{co-planar} pages that are encountered consecutively as one
moves along the arc; Figure \ref{fig7}(c) has a $3$-page book like
this, comprising pages $P_1, P_2$ and $P_3$ and bookcovers $C_1$ and
$C_2$. There may also be a hemi-vertex's book that has co-planarity
of a page and a bookcover encountered on the arc immediately before
or after the page; see Figure \ref{fig2}(c) where page $P$ and
bookcover $C_1$ are coplanar. A non-hemi vertex cannot have a book
with such \emph{consecutive} coplanarity, although non-consecutive
coplanarity is possible. (Figure \ref{fig2}(b) has two books each
with coplanar, non-consecutive, bookcovers.)

{\sc Lemma 2:} \emph{If a vertex $v$ has a $p$-page book associated
with it, then $v$ lies in the interior of $p+1$ ridges of cells and
in the interior of $p$ sides of plates that lie within the book,
with the exception that a hemi-vertex having one of the `consecutive
coplanarity' structures mentioned above is adjacent to $p$ such
ridge-interiors and $p$ such plate--side-interiors. }

{\sc Proof:} (i) Each page of the book is a plate of the \tes\ with
$v \subset$ an interior of a side of the plate and this side
$\subset \ell$. A plate on a bookcover either has a corner that
coincides with $v$ or does not contain $v$. Therefore  $v$ lies in
the interior of $p$ plate--sides within the book. (ii) As one moves
along the arc from one bookcover to the other, one encounters $p+1$
cells of the \tes. If  there is no `consecutive coplanarity', then
each of these cells  has a ridge $\subset \ell$, with $v$ being in
the interior of that ridge. If there is a consecutive coplanarity,
and there can be at most one, the cell involved does not have a
ridge $\subset \ell$ and therefore $v$ is in the interior of just
$p$ ridges.  \hfill$\square$.

{\sc Remark 8:}\emph{ A book without pages is possible, but then the
dihedral angle between the two bookcovers must be $<\pi$, thereby
guaranteeing a ridge along the spine. In such a structure, the
bookcovers are \textit{consecutively} encountered on the arc, of
course. If  we increased the dihedral angle to $\pi$, then this
ridge disappears, the book ceases to exist and the erstwhile
bookcovers (though still considered to be \textit{consecutive})
cease individually to be facets and therefore lose the status of
bookcovers. In summary, it is not possible to have a book with no
pages and two coplanar bookcovers.}

Note that a vertex may have no books associated with it. For
example, none of the vertices in a \ftf\ \tes\ have associated books
and many of our non \ftf\ examples have non-book vertices; every
vertex in the stratum models of Example 9 has no associated book. Of
course, a vertex with no books is not contained in any
ridge-interior  (and vice versa).

Vertices can have more than one book. The total numbers,
$m_{Z_1}\cn(v)$ and $m_{P_1}\cn(v)$, of ridge-interiors and
plate--side-interiors adjacent to $v$ is simply the addition of such
adjacencies over all books associated with $v$. An example is the
schematic of Figure \ref{fig8}(c), where $m_{Z_1}\cn(v) = 3+2+3$ and
$m_{P_1}\cn(v)=2+1+2.$

Multi-book vertices have many variants. For example, Figure
\ref{fig2}(b) shows a ($p=2$)-book above the plane of the horizontal
plates and an obliquely-oriented ($p=1$)-book below with a subtle
sharing of the horizontal plates in the making of book-covers. This
sharing can only occur when the books have different spines. A
hemi-vertex can have multiple-books, though all with a common spine.
For example in Figure \ref{fig7}(c), insert a plate in between pages
$P_2$ and $P_3$ with a corner touching $v$, converting
\includegraphics[width=10mm]{scheme2c.eps} into
\includegraphics[width=10mm]{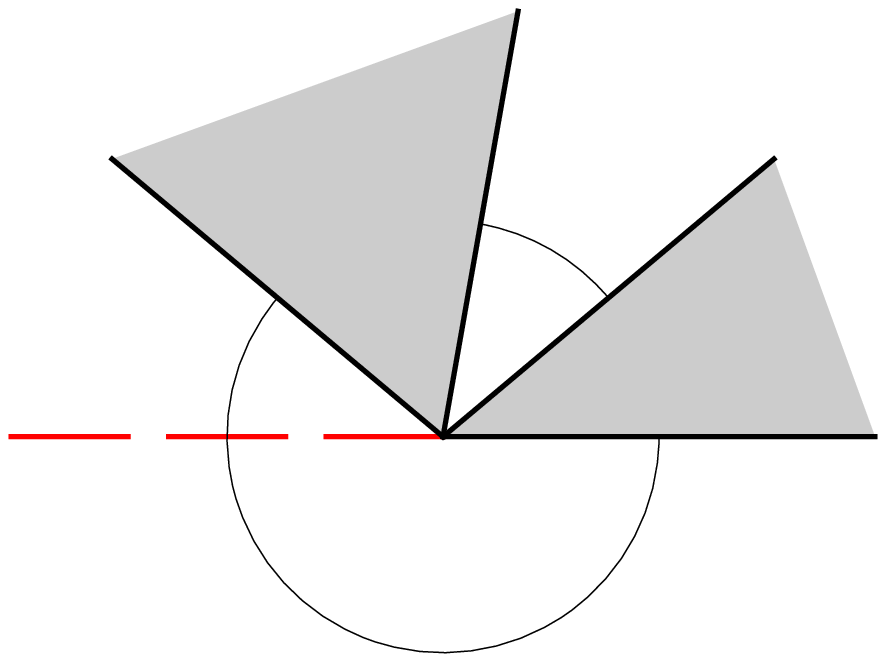}. This creates two books
with common spine, one of the books with no pages and the other with
$P_1$ as its sole page.

In general, the books (if any) for a non-hemi vertex can be mutually
oblique; the generic situation for a non-hemi vertex $v$ with $n$
books is as follows. Draw $n$ lines of differing
orientations passing through $v$. On each line construct a
book initially with no pages using the line as spine. The dihedral angle between the covers will be small if $n$ is large. To increase $p$ for some of these
books, just add further pages.

Figure \ref{fig8}(a) shows an example where $n=2$; the books, each
with $p=0$, are shown with dark shading and one possible version of
the surrounding architecture is being assembled (the fourth
triangular plate almost in position).

Although the additional architecture around $v$ and its books can be very complicated, beyond our powers of visualisation and analysis if $n$ is large and many spines are involved, we know that no contributions to $m_{Z_1}\cn(v)$ and $m_{P_1}\cn(v)$ come from this architecture. The simple additive count over books suffices.

Even when $n=1$, the additional architecture can be complicated. In
Figure \ref{fig7}(c) it is simple, but if that figure is modified by
inserting a plate in between page $P_3$ and book-cover $C_2$ with a
corner touching $v$, it is more elaborate (as shown in Figure
\ref{fig8}(b)). The modification has turned $P_3$ into a book cover
(now labelled $C_3$) which, together with the other book-cover $C_1$
and the pages $P_1$ and $P_2$, form a book (the only book associated
with $v$). The additional architecture comprises the two plates each
with a corner touching $v$ and the structure we previously labelled
$C_2$, although this is no longer a book-cover.

\section*{8. Constraints on the interior parameters $\xi, \kappa, \psi$ and $\tau$}

The considerations in the last section can be applied to achieve
inequalities involving the interior parameters; we do so in Lemma 3
below. Part of this lemma, formula (\ref{xptk}), establishes one of
the bounds for $(\kappa, \xi)$ given  in Theorem 1. Furthermore,
formula (\ref{pt}) provides a component to both the upper bound and
lower bound of inequality (\ref{taubound}) in that theorem.

{\sc Lemma 3:} \emph{For stationary random \tes s of $\R^3$ that
are not \ftf,}
\begin{align}
    0 \leq&\ \psi-\tau \leq \tfrac12 \ \VE \qquad \qquad\mathrm{and}
      \label{pt}\\
    \xi \geq&\ \frac{2(\psi-\tau)+3\kappa}{\VE}= L1 (say),\label{xptk}
\end{align}
\emph{augmenting the basic condition $(\kappa,\xi) \in [0,1]\times (0,1]$
which follows from the definitions. Note: the label $L1$ identifies
this as the first of two lower bounds that we establish for $\xi$.}

{\sc Proof:} We can readily conclude that $\tau \leq \psi$ using
Lemma 2 which shows that, for every book associated with the vertex
$v$, the number of adjacent plate--side-interiors is $\leq$ the
number of adjacent ridge-interiors. From the discussion in Section
7, we have $m_{P_1}\cn(v)\leq m_{Z_1}\cn(v)$ for all $v$. Taking
expectations for typical $v$ shows that $\tau \leq \psi$.

To prove (\ref{xptk}), we use $\xi \mu_{VE} = \mu_{VE [\pi]}$, see (\ref{xipi}), and show that, for a vertex $v\in V$, there are lower bounds on
$m_{E[\pi]}(v)$, the number of $\pi$-edges adjacent to (that is,
emanating from) the vertex $v$.
\begin{equation}\label{p1}
    m_{E[\pi]}(v)\geq
\begin{cases}
    2(m_{Z_1}\cn(v) -m_{P_1}\cn(v))&\qquad \mathrm{if\ }v
       \mathrm{\ is\ not\ a\ hemi-vertex},\\
    4+2(m_{Z_1}\cn(v)-m_{P_1}\cn(v))&\qquad\mathrm{if\ }v
       \mathrm{\ is\ a\ hemi-vertex\ and\ }m_{Z_1}\cn(v)>0,\\
    3&\qquad\mathrm{if\ }v
       \mathrm{\ is\ a\ hemi-vertex\ and\ }m_{Z_1}\cn(v)=0.
\end{cases}
\end{equation}
To show these inequalities, first let $v$ be a non-hemi vertex. With
Lemma 2 the number of books that $v$ has is  $m_{Z_1}\cn(v)
-m_{P_1}\cn(v)$, because none of these books have consecutive
coplanarity. Every book-cover contributes at least one $\pi$-edge to
the count of $ m_{E[\pi]}(v)$. Outside these books there are no
further ridge or plate--side interiors containing $v$. The
arrangement changing Figure \ref{fig7}(c) to Figure \ref{fig8}(b)
(that is, \includegraphics[width=10mm]{scheme2c.eps} to
\includegraphics[width=10mm]{scheme8b.eps}) affects no change of
$m_{Z_1}\cn(v) -m_{P_1}\cn(v)$, whereas the number of $\pi$-edges
adjacent to $v$ can increase; it increases by two in Figure
\ref{fig8}(b), but this could be more with complicated structures
and it might not increase at all (as can be imagined if the two
white plates in Figure \ref{fig8}(b) were aligned). This proves the
first inequality in (\ref{p1}).

Now consider $v$ as a hemi-vertex with one spine and $n\geq1$ books.
If none of these books are \emph{consecutive coplanar}, eliminating
for example

\begin{description}
   \item  \hspace{4cm} \includegraphics[width=18mm]{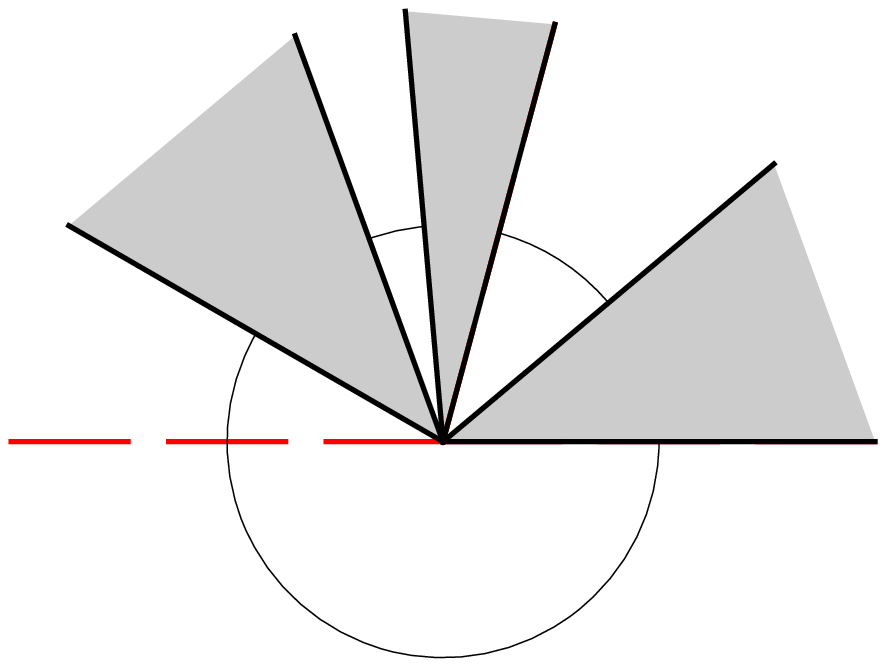} and
\includegraphics[width=18mm]{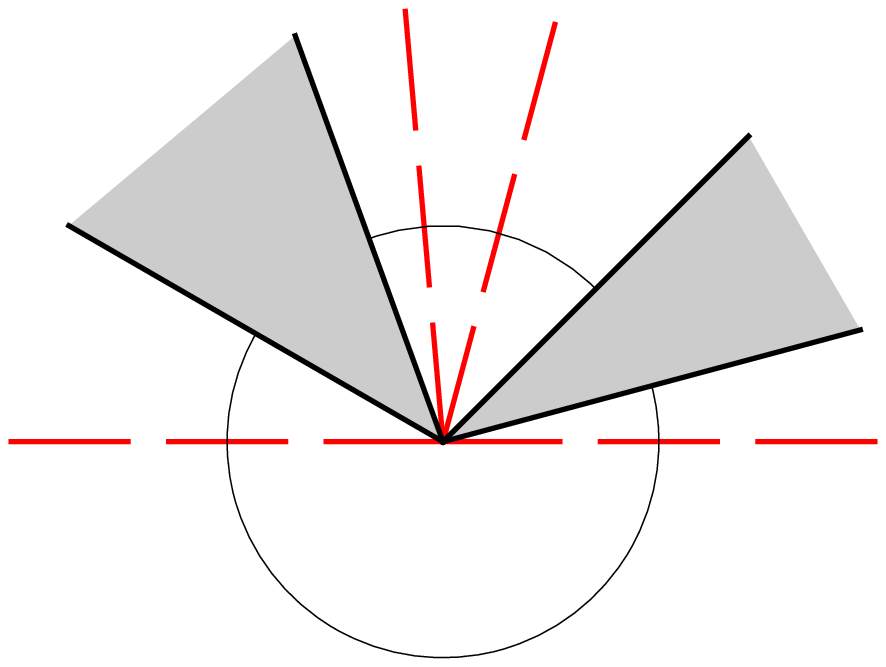} but allowing
\includegraphics[width=18mm]{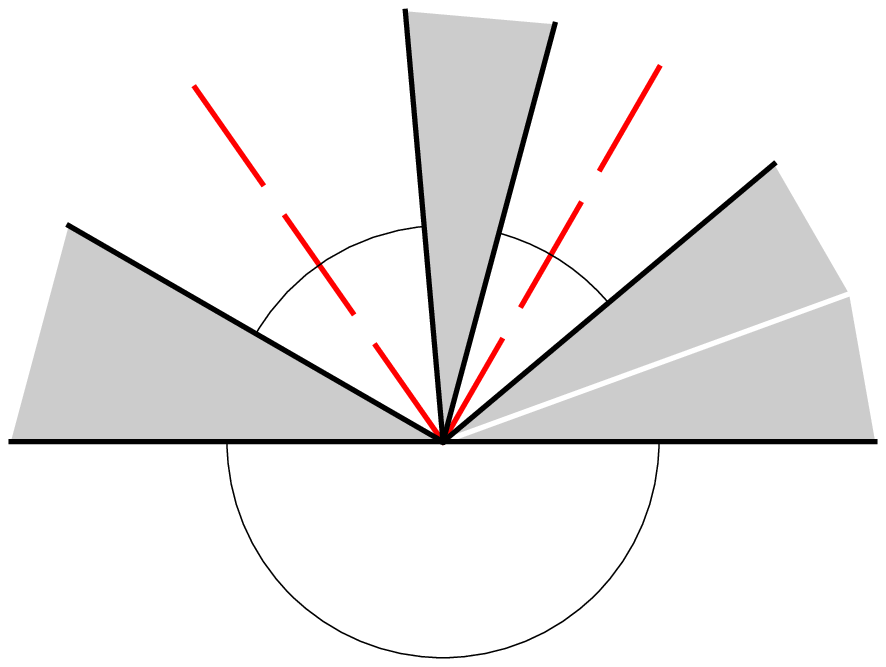},
\end{description}
 the number of
$\pi$-edges emanating from $v$ and contained in the bookcovers is
again $\geq 2(m_{Z_1}\cn(v) -m_{P_1}\cn(v))$. The facet containing
the hemi-vertex $v$  contains at least four \pies\ adjacent to $v$
(two on the spine and at least two off the spine). If one of the
books is consecutive coplanar the number of books now is
$m_{Z_1}\cn(v) -m_{P_1}\cn(v)+1$. The facet containing $v$ has at
least two $\pi$-edges emanating from $v$, see Figure 7(c).

A hemi-vertex with no books must have  three or more emanating $\pi$-edges lying in the interior of the facet that contains $v$. Now (\ref{p1}) is fully proven.

Inequalities  (\ref{p1}) can be written as:
\begin{equation*}
    m_{E[\pi]}(v)\geq
\begin{cases}
    3m_{Z_2}\cn(v)+2(m_{Z_1}\cn(v) -m_{P_1}\cn(v))&\qquad \mathrm{if\ }v
       \mathrm{\ is\ not\ a\ hemi-vertex},\\
       1+3m_{Z_2}\cn(v)+2(m_{Z_1}\cn(v) -m_{P_1}\cn(v))
    &\qquad\mathrm{if\ }v
       \mathrm{\ is\ a\ hemi-vertex\ and\ }m_{Z_1}\cn(v)>0,\\
    3m_{Z_2}\cn(v)+2(m_{Z_1}\cn(v) -m_{P_1}\cn(v))&\qquad\mathrm{if\ }v
       \mathrm{\ is\ a\ hemi-vertex\ and\ }m_{Z_1}\cn(v)=0.
\end{cases}
\end{equation*}
So, for general $v$,
\begin{equation*}
    m_{E[\pi]}(v)\geq 2(m_{Z_1}\cn(v)
    -m_{P_1}\cn(v))+3m_{Z_2}\cn(v).
\end{equation*}
Therefore,
\begin{equation*}
    \mu_{VE[\pi]}\geq 2(\psi
    -\tau)+3\kappa.
\end{equation*}
With (\ref{xipi}) we prove (\ref{xptk}). Inequality (\ref{xptk})
leads to the proof of the upper bound in (\ref{pt}). $\hfill\square$

\textbf{Known results involving $\xi, \kappa, \psi$ and $\tau$:}
Some of the most important results in our earlier paper, \cite{wc},
are three formulae that yield the expected numbers of apices, ridges
and facets for a typical cell --- entities that are equal to the
mean adjacencies $\ZV, \ZE$ and $\ZP$ (given in Table 1) in the
\ftf\ case, but not in general.  These three formulae follow:
\begin{align}
    \nuZA=&\ \ZV - \PV\frac{2(\kappa+\psi)}{f(\PV)};\label{nuZA}\\
    \nuZR =&\ \ZE - \PV\frac{\xi\VE + 2\psi}{f(\PV)};\label{nuZR}\\
    \nuZF =&\ \ZP - \PV\frac{\xi\VE - 2\kappa}{f(\PV)}.\label{nuZF}
\end{align}
Our earlier paper also presented formulae for the mean number of
sides (also of corners) for the typical cell-facet and for the
typical plate.
\begin{align}
    \nuFS = \nuFC =&\ 2\PV \frac{\VE(\EP-\xi)
        -2\psi}{2\VE\EP-\PV(\xi\VE-2\kappa))}\label{nuFS}\\
    \nuPS = \nuPC =&\ \PV\Bigl(1-
    \frac{2\tau}{\VE\EP}\Bigr).\label{nuPS}
\end{align}
Intensities, shown in Table 2 below, were also given there. These
augment the list of intensities for the primitive elements given in
Table 1.
\begin{center}
\begin{tabular}{|c||c|}
\hline $X$ & $\lambda_X/\lambda_V$ \\ \hline \hline
 facets $Z_2$&$\ds{\frac{\overset{\ }{2}\VE\EP-\PV(\xi\VE-2\kappa)}{2\PV}}$\\[3mm]
ridges $Z_1$ & $\frac{1}{2}(\VE(\EP-\xi)-2\psi)$ \\[1mm]
apices $Z_0$ & $\tfrac{1}{2} f(2)- \kappa - \psi$ \\[1mm]\hline
facet--sides $(Z_2)_1$ or facet--corners $(Z_2)_0$&
$\VE(\EP-\xi)-2\psi$
\\[1mm]\hline
plate--sides $P_1$ or plate--corners $P_0$& $\overset{\
}{\frac{1}{2}}(\VE\EP - 2 \tau)$
\\[1mm]\hline
\end{tabular}

\small{Table 2:} {\scriptsize Intensities for non-primitives in
spatial \tes s.  Note the iterated notation, for example $(Z_2)_1$
for facet--sides.  }
\end{center}

We have systematically examined every sensible inequality that can
be applied to these formulae. Firstly, we investigated the
constraints needed on our seven parameters to ensure that all
intensities from Tables 1 and 2 are non--negative (and in some
cases, positive). Secondly we have examined the consequence of the
inequalities $\nuZA\geq 4,\ \nuZR\geq 6, \nuZF \geq 4, \nuFS\geq 3$
and $\nuPS \geq 3$. We do not report details of the constraints from
the first exercise as these were dominated by those from the second,
which yield as follows.
\begin{align}
    \nuZA\geq 4 \quad \Llra&\quad \ZV - \PV\frac{2(\kappa+\psi)}{f(\PV)}\geq
    4\notag\\
    \quad \Llra&\quad \kappa \leq \VE-2+\tfrac12
    \VE\EP\Bigl(1-\frac{4}{\PV}\Bigr)-\psi= K\ (say).\label{K} \\
    \mathrm{Likewise,\ \ } \nuZR\geq 6\quad\Llra&\quad \xi
    \leq 6 (1-\frac{2}{\VE}) +
    \EP\Bigl(1-\frac{6}{\PV}\Bigr)-\frac{2\psi}{\VE}=U1\ (say),\label{U1}\\
    \nuZF\geq 4\quad \Llra& \quad  \xi \leq 4\Bigl(1-\frac{2}{\VE}\Bigr)
    - \frac{2\EP}{\PV}+\frac{2\kappa}{\VE}=U2\ (say)\label{U2}\\
    \mathrm{and\ \ }\nuFS \geq 3\quad \Llra&\quad \xi
    \geq \frac{4\psi+ 6\kappa}{\VE}-2
    \EP\Bigl(1-\frac{3}{\PV}\Bigr)
     = L2\ (say).\label{L2}
\end{align}

In the non \ftf\ theory, the basic identity (\ref{simple})
applicable to each of the \tes's convex cells becomes $2 \nuZR \geq
3 \nuZA$ and this is equivalent to
\begin{equation}\label{U3}
     \xi \leq 3 -\frac{\EP}{2}+ \frac{\psi - 6+3\kappa}{\VE}= U3\ (say),
\end{equation}
an identity which can also be derived from $3\la_Z \nuZA \leq
\la_{Z_2} \nu_0(Z_2)$. The other polyhedral inequality
(\ref{simplex}) yielded $2\nuZR \geq 3 \nuZF$ but this led to
something already known, namely (\ref{L2}).

\textbf{Constraints on $\xi$ and $\kappa$, given $\VE, \EP, \PV,
\psi$ and $\tau$:} We see in (\ref{xptk}) and (\ref{K})--(\ref{U3}) two lower bounds
and three upper bounds for $\xi$, plus one upper bound for $\kappa$.
The bounds are expressed in terms of $\VE, \EP, \PV, \psi$ and
$\tau$ which, for the moment pending further investigation, are
assumed to be `appropriate' --- permitting a non-null set of
$(\kappa, \xi)$ values.

\begin{figure}[ht]
\psfrag{x}{$\xi$}\psfrag{k}{$\kappa$}\psfrag{U1}{\tiny $U1$}
\psfrag{C}{\hspace{-.6mm}\tiny $C$}\psfrag{K}{\hspace{-.6mm}\tiny $K$}
\psfrag{U2}{\hspace{-1mm}\tiny $U2$}
\psfrag{L1andU3}{\hspace{-1.2cm}\tiny $L1$ and $U3 \rightarrow$}
\psfrag{U3}{\hspace{-.7mm}\tiny $U3$}
\psfrag{L1}{\hspace{-.7mm}\tiny $L1$}
\psfrag{L2}{\hspace{-.7mm}\tiny $L2$}
    \begin{center}
    \includegraphics[width=52mm]{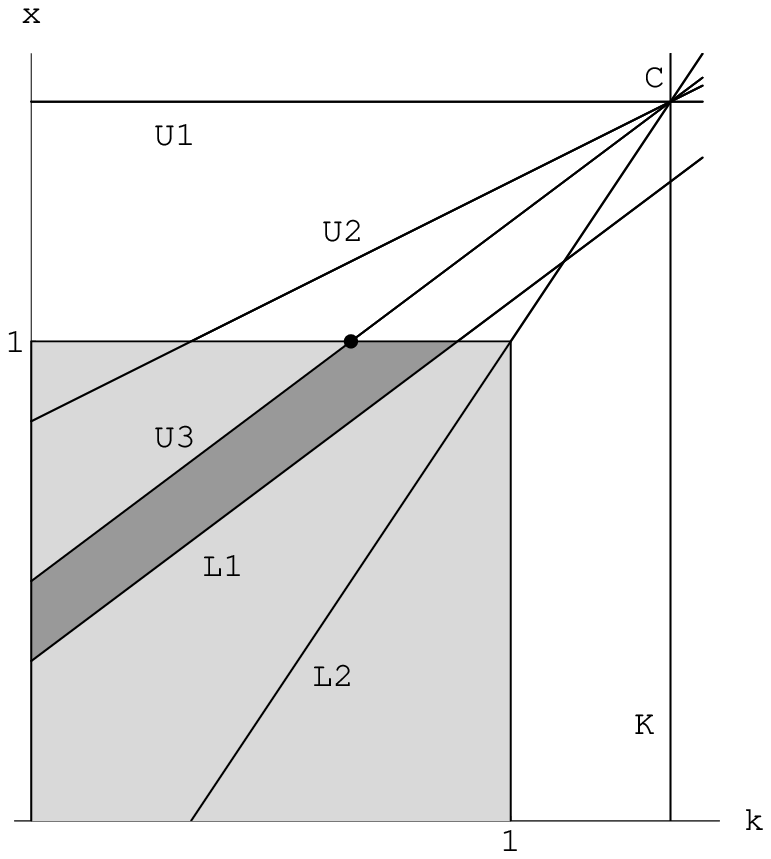}
    \includegraphics[width=52mm]{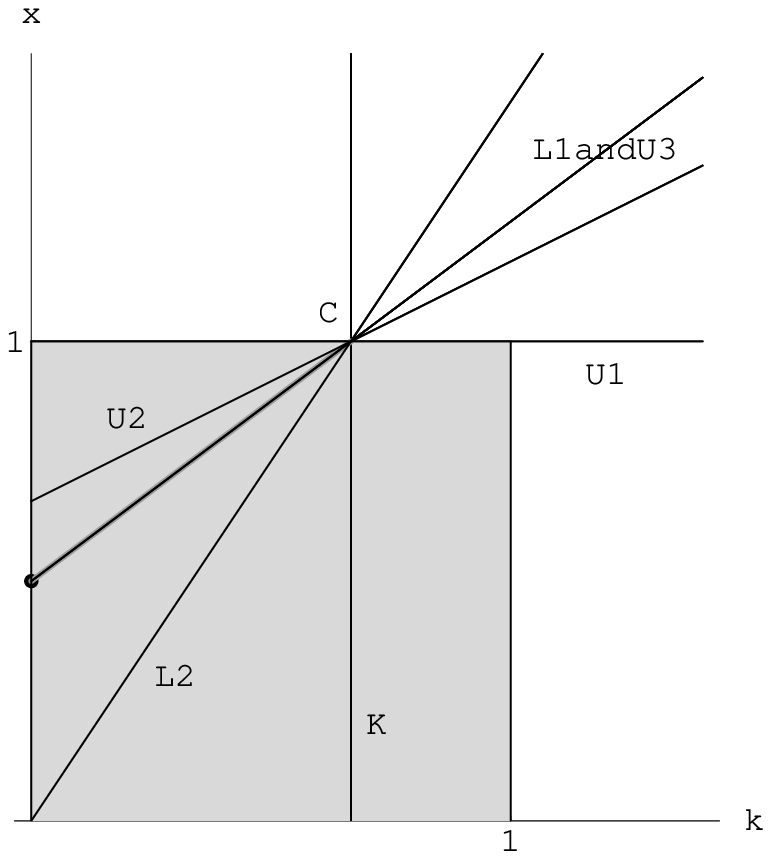}
    \includegraphics[width=48mm]{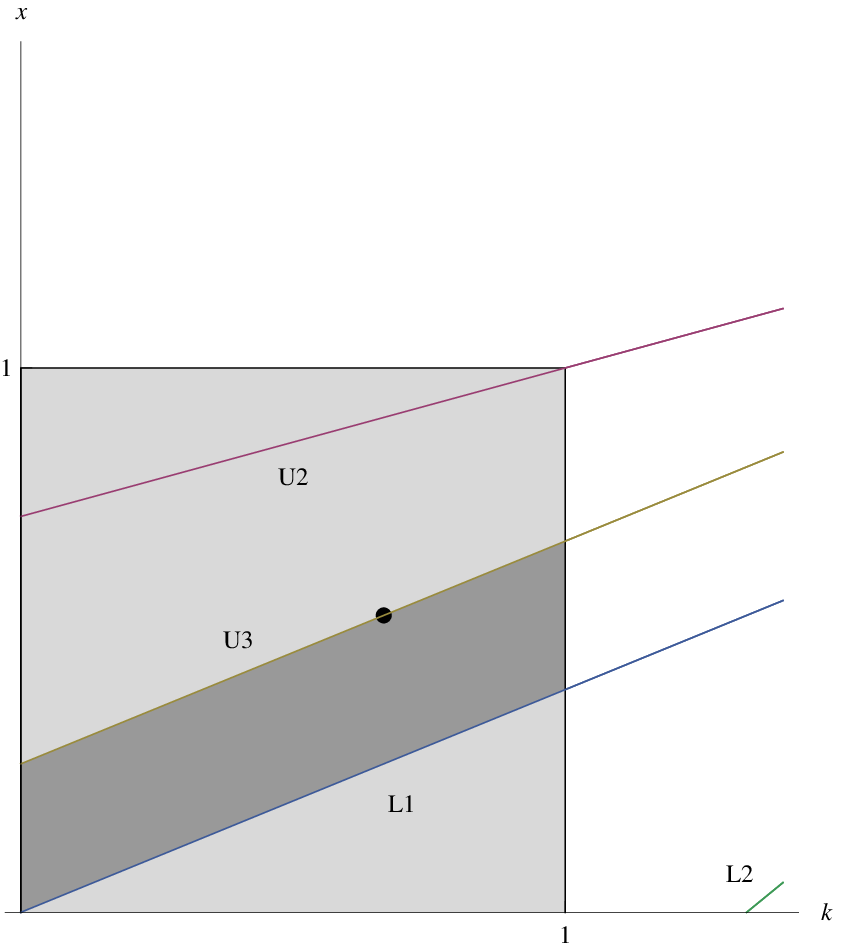}\\
          {\scriptsize (a) Example 4, STIT\hspace{1.5cm}
          (b) Example 6a, triangular prism model\hspace{1.5cm}(c) Example 17}\hspace{2cm}$\ $\\[2mm]
           \includegraphics[width=52mm]{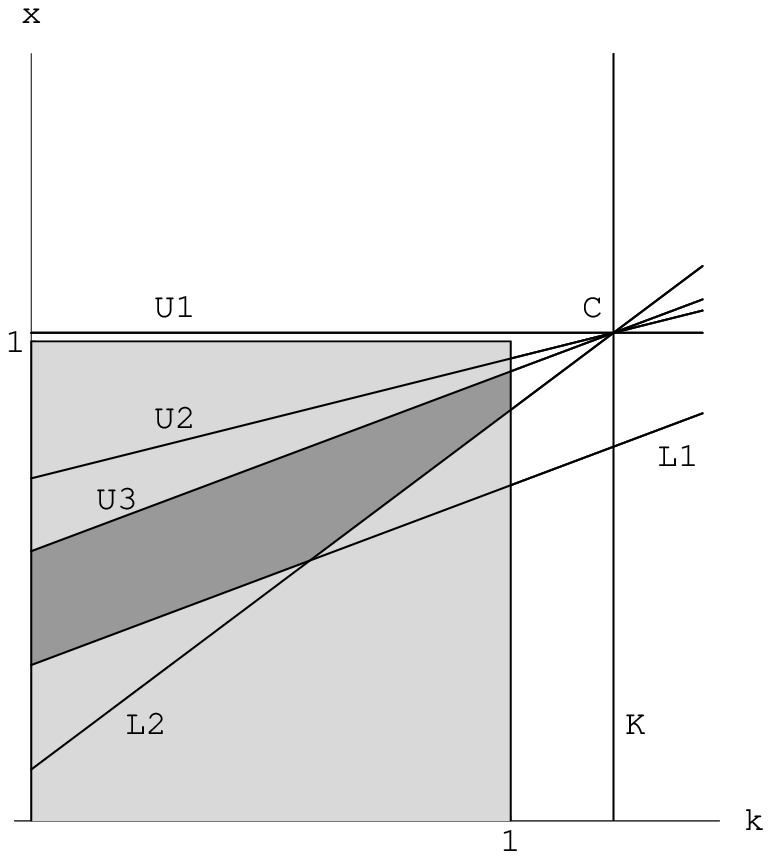}
    \includegraphics[width=52mm]{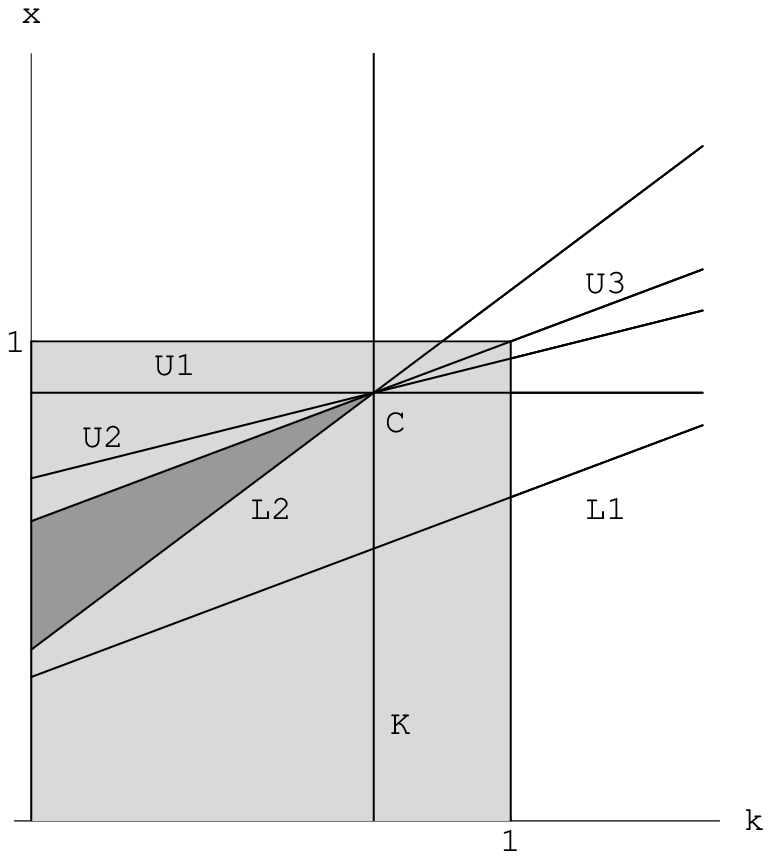}
    \includegraphics[width=52mm]{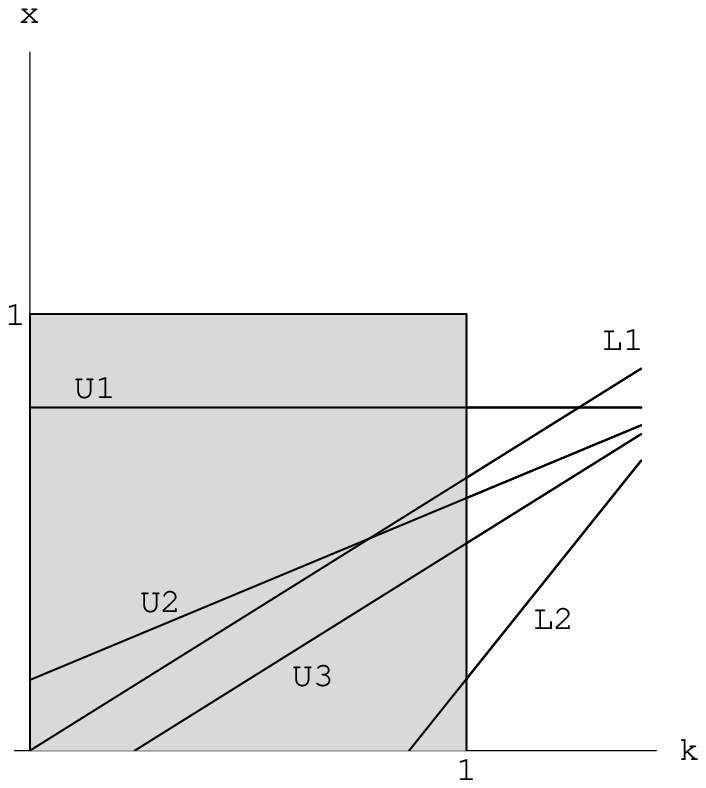}\\
            {\scriptsize (d) $L2$ impacts on zone \hspace{1cm}
          (e) $L1$ has no impact. Also $C \in [0,1]\times (0,1]$  \hspace{1cm}
          (f) Null zone}\hspace{2cm}
        \caption{\label{fig9}\scriptsize  Examples of permitted
        zones for $(\kappa, \xi)$ shown in the darker shading, against
        the light grey region $[0,1]\times(0,1]$. In (a)-(c), three of our examples
        are featured, whilst (d)-(f) show other variations. For (a)-(c), the heavy dot shows the $(\kappa, \xi)$ point of that example whilst the dark shading shows the permitted range of $(\kappa, \xi)$ for \tes s which have the same $\VE, \EP, \PV$ and $\psi$ values as the model which produced the heavy dot.
        In (b), the darker zone is
        a line--segment running between the heavy dot and $C$. In (f) the zone is null because the lower bound
         $L1$ is higher
        than the upper bound $U3$, implying that the choice of $\VE, \EP, \PV, \psi$ and
$\tau$ is inappropriate. For the record, (d) and (e) use $\VE=8,
\EP=4$ and $\PV=\tfrac72$ (with $(\psi, \tau) = (\tfrac52,\tfrac65)$
for (d) and $(\psi, \tau)= (3, \tfrac95)$ for (e)). In (f),
$\VE=\tfrac{24}{5}, \EP=\tfrac{19}{5}$ and  $\PV=\tfrac72$, with
$\psi=\tau=0$.}
    \end{center}
\end{figure}

Considerable insight on the permitted $(\kappa, \xi)$ values comes
from Figure \ref{fig9}, where the bounds appear as straight lines
--- because the bounds on $\xi$  are linear in $\kappa$.

{\sc Lemma 4:} \emph{The permitted range for $\kappa$ and $\xi$,
given $\VE, \EP, \PV, \psi$ and $\tau$, is given by $0\leq\kappa\leq
\min[1,K],\, \xi >0$ and $\max(L1,L2)\leq \xi \leq \min[1,U3]$,
written in full form as}
\begin{align*}
        0\leq \kappa \leq\ \min\Bigl[1,\VE-2+\tfrac12
  \VE\EP
  \Bigl(1-\frac{4}{\PV}\Bigr)&-\psi\Bigr]
  \qquad\mathrm{and}\\
    \max\Bigl[\frac{2(\psi-\tau)+3\kappa}{\VE},\frac{4\psi+ 6\kappa}{\VE}-2
    \EP\Bigl(1-\frac{3}{\PV}\Bigr)\Bigr] \leq&\ \xi\ \leq \min\Bigl[1,
    \ 3 -\frac{\EP}{2}+ \frac{\psi - 6+3\kappa}{\VE}\Bigr],
\end{align*}
\emph{supplemented by} $\xi >0$.

{\sc Proof:} Figure \ref{fig9} suggests that five of the straight
lines pass through a common point, labelled $C$, whose coordinates
are obviously $(\kappa,\xi) =(K,U1)$. It can be proved generally
with simple algebra that $C$ lies on the other three lines
 $U2, U3$ and $L2$ (details omitted). The redundancy
of $U1$ and $U2$ follows, because the line for $U3$ has a higher
slope than those for $U_1$ and $U_2$, implying $U_1\geq U_3$ and
$U_2\geq U_3$ when $\kappa \leq K$ (as it always is, from
(\ref{K})).

Bound $L1$ is established in Lemma 3. The redundancy of $U1$ and $U2$, leaves only (\ref{K}),(\ref{L2}) and (\ref{U3}), and they provide the remaining bounds to complete the lemma. \hfill$\square$.

Note that the permitted zone for $(\kappa,\xi)$ can have zero area,
as in Figure \ref{fig9}(b) where the zone is a line--segment or in
Example 8, the Divided Delaunay \tes, where one can show that the
zone is the single point $(\kappa, \xi) = (0, 64\pi^2/(35+112\pi^2))
\approx (0, 0.554)$.

Figure \ref{fig9}(f) raises an interesting issue. The lines $L1$ and
$U3$ can never cross because they are parallel. Yet, the values
chosen for $\VE, \EP, \PV, \psi$ and $\tau$ have created the
situation where $L1>U3$. This is clearly not allowed, so something
must be inappropriate in the choice of $\VE, \EP, \PV, \psi$ and
$\tau$. This issue motivates our next subsection.

\textbf{Appropriate values for $(\psi,\tau)$:} We must avoid a
choice of $\psi$ and $\tau$ which leads to $L1>U3$. For a \tes\ to
exist, we require $L1\leq U3$ which after rearrangement becomes
\begin{equation}\label{tau}
    \tau \geq \frac{\psi}{2} +\frac{\VE}{4}\Bigl(\EP-6\bigl(1-\frac{2}
    {\VE}\bigr)\Bigr).
\end{equation}
This is the first example of a constraint on $(\psi,\tau)$ arising from the need to avoid a null space for $(\kappa, \xi)$. Other such constraints follow shortly, but firstly we consider
one other inequality which  applies for $\tau$. From inequality (\ref{nuPS}), which represents the requirement $\nuPS\geq 3$, we have that
\begin{equation}\label{tau1}
\tau \leq \tfrac12\VE\EP\Bigl(1-\frac{3}{\PV}\Bigr).
\end{equation}
The bound (\ref{tau1}) is the last inequality (used in conjunction with (\ref{pt}) of Lemma 3) needed to complete our proof of (\ref{taubound}) in Theorem 1.

{\sc Remark 9:} \emph{Consider again (\ref{tau}). Note that the term
in the larger brackets relates to the fundamental curve of Figure
\ref{fig1}. So, the \tes\ is represented by a point above the curve
in Figure \ref{fig1}, if and only if $\tau>0$. Thus the only \tes s
below the curve have $\tau=0$.  These findings were anticipated in
the caption of Figure \ref{fig5}.}

There does not exist a \tes\ if $K<0$ or if $U1<0$ and this fact
constrains $\psi$ somewhat. From (\ref{K}) and (\ref{U1}), we have
respectively two constraints:
\begin{align}
    \psi \leq&\ \VE-2+\frac{\VE\EP}{2}\Bigl(1-\frac{4}{\PV}\Bigr);\label{psi1}\\
    \psi \leq&\ 3\VE-6+\frac{\VE\EP}{2}\Bigl(1-\frac{6}{\PV}\Bigr)\label{psi2}.
\end{align}
Their violation would make the $(\kappa, \xi)$ domain null, so they are needed (or at least (\ref{psi1}) is).

{\sc Lemma 5:} \emph{For non \ftf\ \tes s, the bound in (\ref{psi2}) is greater than that in (\ref{psi1}) and is therefore redundant.  }

{\sc Proof:} From (\ref{ZV4}), we know that $\PV > \tf12 \frac{\VE\EP}{\VE-2}$ when the \tes\ is not \ftf. Simple algebra shows that  $\PV > \tf12 \frac{\VE\EP}{\VE-2} \Llra$ ``bound of (\ref{psi2}) $>$ bound of (\ref{psi1})".  \hfill$\square$

We can also derive another upper bound for $\psi$ from (\ref{L2}). Although $L2>U3$ is not possible for any $\kappa$, so no tessellation--existence question arises in this way, there would be an issue if $\{\xi\geq L2\}\cap \{(\kappa, \xi) \in [0,1]\times (0,1]\} =\emptyset$. In view of $L2$'s positive slope, it is enough to require that $L2\leq 1$ when $\kappa=0$ to avoid this problem and this yields the inequality
\begin{align}
    \psi\leq
      &\ \frac{\VE}{4}+\frac{\VE\EP}{2}\Bigl(1-\frac{3}{\PV}\Bigr).\label{psi3}
\end{align}
No further nullity-avoidance constraints are required, because Remark 10 below eliminates such issues with $L1$ and $U3$.

{\sc Remark 10:} \emph{With $0  \leq  L1 \leq 1$ when $\kappa=0$,
the intersection of  $\{(\kappa,\xi)\in [0,1]\times (0,1]\}$ and
$\{\xi\geq L1\}$ is non-empty. See six cases of $L1$ in Figure
\ref{fig9} which illustrate the line $L1$. If $\psi -\tau = \tf12
\VE$, then the only permitted case is $(\kappa,\xi) = (0,1)$.}

Simple algebraic manipulation establishes the following lemma which in turn establishes (\ref{psibound}), thus completing the proof for another part of Theorem 1.

{\sc Lemma 6:} \emph{Denote the right-hand sides of
(\ref{psi1}) and(\ref{psi3}) by $R1$ and $R2$. Then}
\begin{align*}
    R2 < R1  \quad \Llra&\qquad \frac{2\VE\EP}{3\VE-8} < \PV \qquad \mathrm{and\ }\\
    R1 \leq R2 \quad \Llra&\quad  \PV \leq \frac{2\VE\EP}{3\VE-8}.
\end{align*}
\emph{Moreover, $\ds{\frac{2\VE\EP}{3\VE-8}}$ lies within the permitted range for $\PV$ shown in (\ref{case2}).}

We conclude this subsection having established five constraints on the $(\psi, \tau)$ domain that have not been made redundant, namely $R1$ and $R2$ of Lemma 6 and (\ref{pt}), (\ref{tau}) and (\ref{tau1}).

\begin{figure}[ht]
    \psfrag{EP}{{ \scriptsize$\hspace{-4mm}\EP$}}
        \psfrag{VE}{{ \scriptsize$\hspace{-4mm}\VE$}}
    \psfrag{PV}{{ \scriptsize$\hspace{-3.4mm}\PV$}}
    \psfrag{tau}{$\tau$} \psfrag{psi}{\hspace{-2mm}$\psi$}
    \psfrag{Ro}{\scriptsize$R1$}  \psfrag{Rt}{\scriptsize$R2$}
    \psfrag{eq}{\scriptsize$\hspace{-2mm}\tau = \psi$}
    \psfrag{17}{\scriptsize\hspace{-3mm} (\ref{pt})}
        \psfrag{30}{\scriptsize \hspace{-3mm}(\ref{tau})}
        \psfrag{22}{\scriptsize (\ref{tau1})}
    \begin{center}
    \includegraphics[width=170mm]{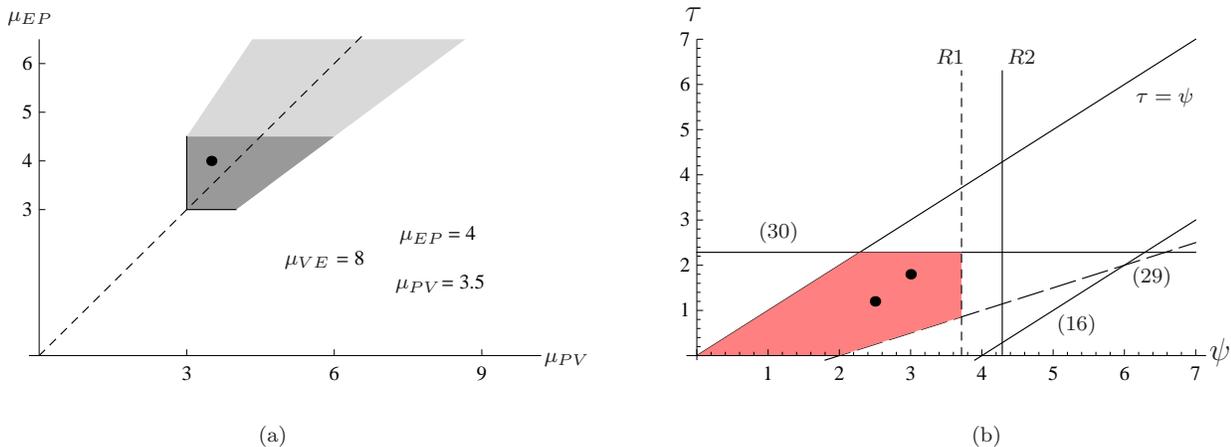}\\
          {\scriptsize (a) \hspace{9cm}(b)  }
        \caption{\scriptsize  \label{fig10} (a) This follows the style in Figure \ref{fig5}, except now we add the dashed line, whose equation is $\EP = \tf12(3\VE-8)\PV/\VE$. The dot corresponds to $(\PV, \EP) = (\tf72, 4))$.  (b) The permitted domain for $(\psi, \tau)$ given $(\PV, \EP) = (\tf72, 4))$, with the various bounds marked. See the text for further explanation and the connection with Figure \ref{fig9}. }
    \end{center}
\end{figure}

\textbf{The permitted $(\psi, \tau)$ domain; an example plot:} In
parts (d) and (e) of Figure \ref{fig9} are two examples sharing
common values for the cyclic parameters, namely $\VE=8,\, \EP=4$ and
$\PV=\tfrac72$. Figure \ref{fig10}(a) shows the allowed domain for
$(\PV, \EP )$, given $\VE=8$ (calculated from Theorem 1), together
with the dot corresponding to the particular case. This drawing is
in essentially the same format as the plots in Figure \ref{fig5}.
The dashed line is an innovation, separating the domain into two
parts, in keeping with Lemma 6 and (\ref{psibound}). The particular
dot lies to the left of this dashed line, indicating that $R1<R2$
and this is confirmed in Figure \ref{fig10}(b), where $R1$ and $R2$
are marked.

In  Figure \ref{fig10}(b), we show the allowed $(\psi, \tau)$ domain
given the particular values of all three cyclic parameters. We see
the all the constraints on $\tau$ and $\psi$ developed earlier in
this section. The two dots in this domain are the values used in the
preparation of Figure \ref{fig9}, with $(\psi, \tau) =
(\tfrac52,\tfrac65)$ for (d) and $(\psi, \tau)= (3, \tfrac95)$ for
(e). Note that each dot produces a qualitatively different $(\kappa,
\xi)$ domain as seen in Figure \ref{fig9}.

\textbf{Are there further constraints on the cyclic parameters?} We have not yet addressed the potential for the $(\psi, \tau)$ domain to be empty. Perhaps some values of $\VE, \EP$ and $\PV$ which we have to date regarded as valid, create this null situation and lead to the non-existence of a \tes? These concerns can, however, be dismissed.

{\sc Lemma 7: } \textit{Cyclic parameters satisfying (\ref{case2}) in  Theorem 1  cannot lead to a null domain for $(\psi, \tau)$. }

{\sc Proof:} Simple algebraic calculations show that $\tf12 \VE\EP
(1-3/\PV)$, which is the right-hand side of (\ref{tau1}), is less
than $\min(R1,R2)$. It is also easy to show that $\tf12 \VE\EP
(1-3/\PV)$ is greater than the maximum of the two lower bounds for
$\tau$ given in (\ref{pt}) and (\ref{tau})
 when both these bounds are evaluated at $\psi=
\min(R1,R2)$. It is also obvious that $\min(R1,R2)>0$. The
statements above hold for all values of the cyclic parameters which
satisfy (\ref{case2}). This proves that the domain for $(\psi,
\tau)$ is not null under the stated premise of the lemma.
\hfill$\square$

\begin{figure}[ht]
\psfrag{tau}{$\tau$} \psfrag{psi}{\hspace{-2mm}$\psi$}
    \psfrag{EP}{{ \scriptsize$\hspace{-4mm}\EP$}}
        \psfrag{VE}{{ \scriptsize$\hspace{-4mm}\VE$}}
    \psfrag{PV}{{ \scriptsize$\hspace{-3.4mm}\PV$}}
    \begin{center}
    \includegraphics[width=180mm]{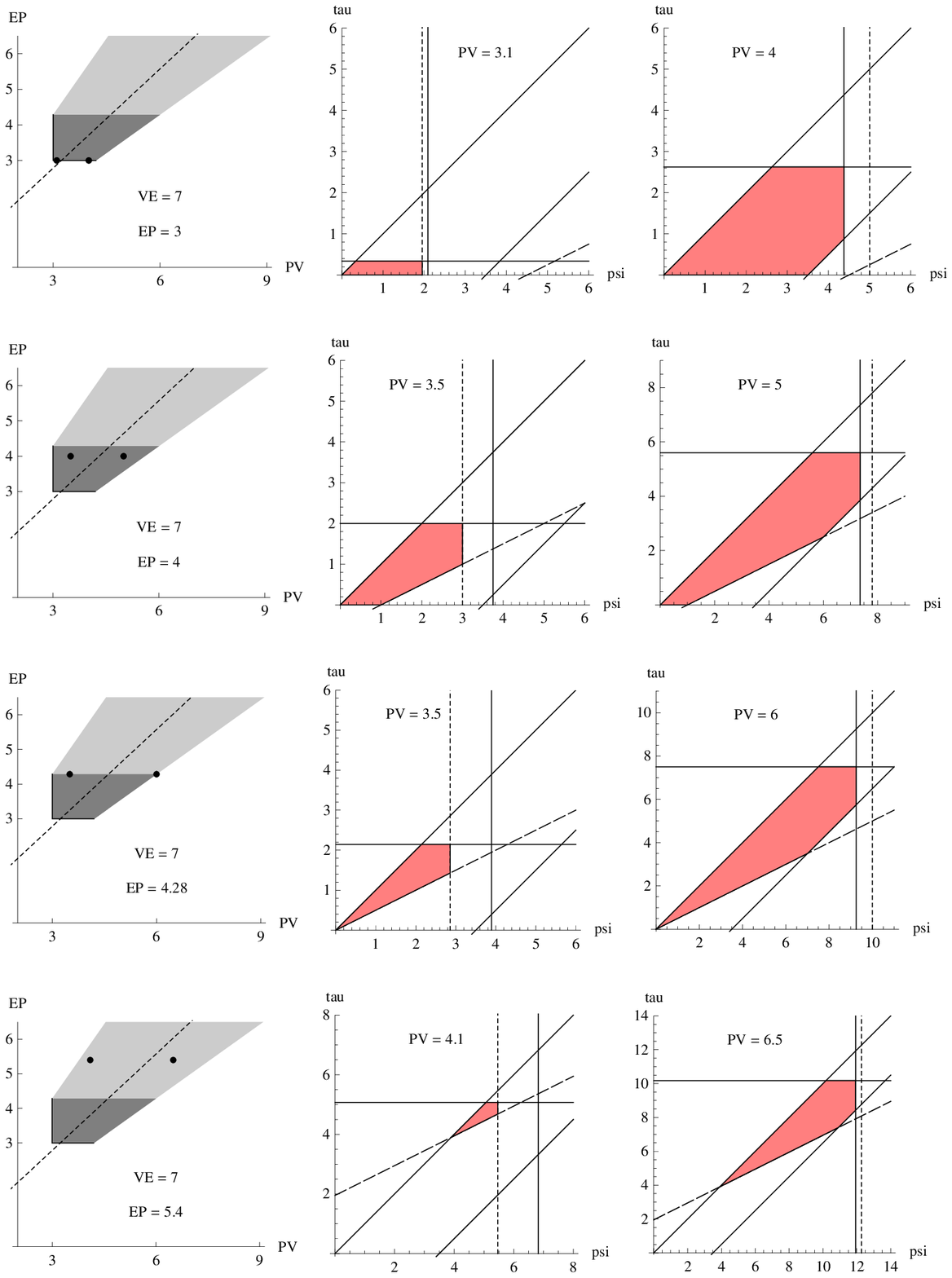}\\
        \caption{\scriptsize  \label{fig11} A diversity of shapes for the $(\psi,\tau)$ permitted domain when $\VE=7$. The dots in the leftmost diagram of each row set the cyclic parameters for the two diagrams to its right. }
    \end{center}
\end{figure}

We conclude this section with the twelve diagrams that comprise
Figure \ref{fig11}; in all diagrams $\VE=7$. Also $\EP$ increases as
one goes down the page. Two values of $\PV$ are shown on each row. A
diversity of shapes for the $(\psi,\tau)$ permitted domain is
demonstrated.

\section*{9. Concluding remarks}

Although it is not easy to work in a seven-dimensional space, we
think that our approach has produced an intelligible picture of the
constraints that apply for spatial \tes s, in both the \ftf\ case
and the general case.

We emphasise however that our focus has been \emph{combinatorial
topology} and our constraints are about \emph{topological}
parameters. We anticipate that further constraints will arise when
geometric aspects are studied in greater detail.

Our approach to the topic of this paper has been rooted in the
literature of \emph{random stationary \tes s} of $\R^3$ and,
although we allude to the \emph{tiling} literature and present many
examples which would be called `tilings', we have not asked the same
questions that tilers ask. We notice, however, that non \ftf\
tilings of $\R^3$  have not often been studied in the literature, so
we think that our findings have value for tilers.

We hope in later work to develop an interface between the two literatures. Tentative steps in that direction have shown that our constraints are useful in showing that tilings by certain polyhedra or combinations of polyhedra do not exist.

At a very late stage of our study, we found three more examples
which allayed a concern that had developed. Up to Example 17, we had
placed no points in our Figure \ref{fig1} on the line $\EP=3$, other
than three having $\VE=4$. As this line is a boundary of our space,
we sought other examples. We found them by using the well-known
monohedral tiling of space  by the rhombic dodecahedron $D$ (with
each facet being a rhombus whose long diagonal is $\sqrt 2$ times
its short diagonal); see \cite{wel}. Using the classical \ftf\
tiling using $D$, we obtained $(\VE, \EP) = (\tf{16}{3},3)$, but
improved on this by cutting $D$ into smaller convex polyhedra.
Importantly these cuts did not hit a ridge-interior of $D$. Two
other non \ftf\ tilings resulted with $(\VE, \EP) = (8,3)$ and
$(\VE, \EP) = (10,3)$. All three are plotted on Figure \ref{fig1} as
18a, 18b and 18c.

{\sc Acknowledgement:} The second author was supported by the German
research foundation (DFG), grant WE 1799/3-1. The first author's
accommodation expenses, while visiting Jena for a period of joint
research, came from the same grant.


\medskip
\begin{thebibliography}{99}

    \bibitem{cow78}
{\sc Cowan, R. (1978).}  The use of ergodic theorems in random
geometry. \it Suppl. Adv. Appl. Prob. \bf 10, \rm 47--57.

    \bibitem{cow80}
{\sc Cowan, R. (1980).}  Properties of ergodic random mosaic
processes. \it Math. Nachr. \bf 97, \rm 89--102.

    \bibitem{cw1}
{\sc Cowan, R. and Weiss, V. (2013).}  Graphical presentations of
the $7$--dimensional parameter space arising in \tes s of $\R^3$.
\emph{Technical note} available as document 85(a) on
\emph{www-personal.usyd.edu.au/$\sim$rcowan/professional/randomgeom.html}

    \bibitem{km}
{\sc Kendall, W. S. and Mecke, J.} (1987). The range of mean-value
quantities of planar \tes s. \emph{J. Appl. Prob.}, \textbf{24},
411--421.

    \bibitem{lz}
{\sc Leistritz, L. and Z\"{a}hle, M.} (1992). Topological Mean Value
Relationships for Random Cell Complexes. \emph{Math. Nachr.}
\textbf{155}, 57--72.

    \bibitem{mec84}
{\sc Mecke, J.} (1984).  Parametric representation of mean values
for stationary random mosaics. \it Math. Operations. Statist. Ser.
Statist. \bf 15, \rm 437--442.

    \bibitem{mil}
{\sc Miles, R. E.} (1971). Poisson flats in Euclidean spaces. Part
II: Homogeneous Poisson flats and the complementary theorem.
\emph{Adv. Appl. Prob.}, \textbf{3}, 1--43.

    \bibitem{mol}
{\sc M{\o}ller, J.} (1989). Random tessellations in $\R^d$. {\em
Adv. Appl. Prob.\/} {\bf 21}, 37--73.

    \bibitem{nw1}
{\sc Nagel, W. and Weiss, V.} (2008). Mean values for homogeneous
STIT tessellation in 3D. \emph{Image Anal. Stereol.}, \textbf{27},
29--37.

    \bibitem{nwc}
{\sc Nguyen, N. L., Weiss, V.  and Cowan, R.} (2013).  Spatial \tes s derived from planar \tes s. In preparation.


    \bibitem{okabe}
{\sc Okabe, A., Boots, B., Sugihara, K. and Chiu, S. N.} (2000).
\emph{Spatial Tessellations: Concepts and Applications of Voronoi
Diagrams}. 2nd. ed., Wiley, Chichester.

    \bibitem{rad}
{\sc Radecke, W.} (1980). Some mean-value relations on stationary
random mosaics in the space. \emph{Math. Nachr.} \textbf{97},
203--210.

    \bibitem{schw}
{\sc Schneider, R. and Weil, W.} (2008). \emph{Stochastic and
Integral Geometry}. Springer, Berlin Heidelberg.


    \bibitem{thw}
{\sc Th\"{a}le, C. and Weiss, V.} (2010). New mean values for
homogeneous spatial \tes s that are stable under iteration.
\emph{Image Anal. Stereol.} \textbf{29}, 143--157.

    \bibitem{wc}
{\sc Weiss, V. and Cowan, R.} (2011). Topological relationships in
spatial \tes s. {\em Adv. Appl. Prob.} {\bf 43}, 963--984.


    \bibitem{wz}
{\sc Weiss, V. and Z\"{a}hle, M.} (1988). Geometric Measures for
Random Curved Mosiacs of $\R^d$. \emph{Math. Nachr.} \textbf{138},
313--326.

    \bibitem{wel}
{\sc Wells, D.} (1991). \emph{The Penguin Dictionary of Curious and
Interesting Geometry}. Penguin Books, London.


    \bibitem{zie}
{\sc Ziegler, G. M.} (2003). Face Numbers of $4$-Polytopes and $3$-Spheres. arXiv:math/0208073v2[math.MG].


\end{thebibliography}
\end{document}